\documentclass[11pt,a4paper]{article}
\pdfoutput=1

\RequirePackage[english]{babel}
\RequirePackage{dutchcal}

\RequirePackage{polski}
\RequirePackage[utf8]{inputenc}
\RequirePackage[T1]{fontenc}

\RequirePackage{amsmath}
\RequirePackage{amssymb}
\RequirePackage{amsthm}
\RequirePackage{mathabx}
\RequirePackage{tikz}
\RequirePackage{tikz-cd}
\RequirePackage{todonotes}
\RequirePackage{hyperref}
\RequirePackage[capitalise]{cleveref}
\RequirePackage{footmisc}
\RequirePackage{bbm}
\RequirePackage{afterpage}
\RequirePackage{mathrsfs}
\RequirePackage{sseq}
\RequirePackage{extarrows}
\RequirePackage{enumitem}
\RequirePackage[heavycircles]{stmaryrd}
\RequirePackage{pdflscape}
\RequirePackage{graphicx}

\newcommand*{\Const}{\mathrm{Const}}
\newcommand*{\Dec}{\mathrm{D\acute{e}c}}
\newcommand*{\Ell}{\mathrm{Ell}}

\newcommand*{\et}{\mathrm{\acute{e}t}}

\newcommand*{\syn}{\mathrm{syn}}

\renewcommand*{\top}{\mathrm{top}}

\DeclareMathOperator{\Spec}{Spec}
\DeclareMathOperator{\Spf}{Spf}

\DeclareMathOperator{\SL}{SL}

\newcommand*{\ASS}{\mathrm{ASS}}

\newcommand*{\DSS}{\mathrm{DSS}}

\newcommand*{\KU}{\mathrm{KU}}
\newcommand*{\KO}{\mathrm{KO}}
\newcommand*{\ko}{\mathrm{ko}}

\newcommand*{\SMF}{\mathrm{SMF}}
\newcommand*{\Smf}{\mathrm{Smf}}
\newcommand*{\smf}{\mathrm{smf}}
\newcommand*{\MU}{\mathrm{MU}}

\newcommand*{\Sph}{\mathbf{S}}

\newcommand*{\TMF}{\mathrm{TMF}}
\newcommand*{\tmf}{\mathrm{tmf}}
\newcommand*{\Tmf}{\mathrm{Tmf}}

\newcommand*{\Aff}{\mathrm{Aff}}
\newcommand*{\CAlg}{\mathrm{CAlg}}

\newcommand*{\DM}{\mathrm{DM}}
\newcommand*{\Fun}{\mathrm{Fun}}
\newcommand*{\Fil}{\mathrm{Fil}}
\newcommand*{\FG}{\mathrm{FG}}

\newcommand*{\Mod}{\mathrm{Mod}}

\newcommand*{\QCoh}{\mathrm{QCoh}}

\newcommand*{\Shv}{\mathrm{Shv}}

\newcommand*{\Stable}{\mathrm{Stable}}
\newcommand*{\Sp}{\mathrm{Sp}}

\newcommand*{\Syn}{\mathrm{Syn}}

\DeclareMathOperator*{\colim}{colim}

\DeclareMathOperator{\map}{map}

\newcommand*{\op}{\mathrm{op}}
\newcommand*{\WH}{\mathrm{Wh}}
\newcommand*{\Wh}{\mathrm{Wh}}

\newcommand{\C}{\mathcal{C}}
\newcommand{\CP}{\mathbf{CP}}
\newcommand{\E}{\mathbf{E}}
\newcommand{\EO}{\mathrm{EO}}

\newcommand{\F}{\mathcal{F}}
\newcommand{\FF}{\mathcal{F}}
\newcommand{\G}{\mathbf{G}}

\newcommand{\M}{\mathsf{M}}

\renewcommand{\O}{\mathcal{O}}

\newcommand{\X}{\mathsf{X}}
\newcommand{\x}{\mathsf{X}}

\newcommand{\y}{\mathsf{Y}}
\newcommand{\Z}{\mathbf{Z}}

\newcommand{\al}{\alpha}

\newcommand{\Ga}{\Gamma}

\DeclareFontFamily{U}{dmjhira}{}
\DeclareFontShape{U}{dmjhira}{m}{n}{ <-> dmjhira }{}

\RequirePackage{mathtools}

\theoremstyle{theorem}\numberwithin{equation}{section}
\newtheorem{theorem}[equation]{Theorem}
\crefname{theorem}{{Th}.\!\!}{{Ths}.\!\!}
\Crefname{theorem}{{Th}.\!\!}{{Ths}.\!\!}
\newtheorem{theoremalph}{Theorem}

\crefname{theoremalph}{{Th}.\!\!}{{Ths}.\!\!}
\Crefname{theoremalph}{{Th}.\!\!}{{Ths}.\!\!}

\Crefname{problem}{{Prb}.\!\!}{{Prbs}.\!\!}
\newtheorem{prop}[equation]{Proposition}
\Crefname{prop}{{Pr}.\!\!}{{Prs}.\!\!}
\newtheorem{lemma}[equation]{Lemma}
\Crefname{lemma}{{Lm}.\!\!}{{Lms}.\!\!}
\newtheorem{cor}[equation]{Corollary}
\Crefname{cor}{{Cor}.\!\!}{{Cors}.\!\!}

\Crefname{conjecture}{{Conj}.\!\!}{{Conjs}.\!\!}

\theoremstyle{definition}\numberwithin{equation}{section}
\newtheorem{mydef}[equation]{Definition}
\Crefname{mydef}{{Df}.\!\!}{{Dfs}.\!\!}

\Crefname{recall}{{Rcl}.\!\!}{{Rcls}.\!\!}
\newtheorem{construction}[equation]{Construction}
\Crefname{construction}{{Con}.\!\!}{{Cons}.\!\!}

\Crefname{ass}{{As}.\!\!}{{As}.\!\!}
\newtheorem{notation}[equation]{Notation}
\Crefname{notation}{{Nt}.\!\!}{{Nts}.\!\!}

\Crefname{situation}{{St}.\!\!}{{Sts}.\!\!}
\newtheorem{variant}[equation]{Variant}
\Crefname{variant}{{Var}.\!\!}{{Vars}.\!\!}

\theoremstyle{remark}\numberwithin{equation}{section}
\newtheorem{example}[equation]{Example}
\Crefname{example}{{Ex}.\!\!}{{Exs}.\!\!}

\Crefname{nonexample}{{NonEx}.\!\!}{{NonEx}.\!\!}

\Crefname{claim}{{Clm}.\!\!}{{Clms}.\!\!}
\newtheorem{remark}[equation]{Remark}
\Crefname{remark}{{Rmk}.\!\!}{{Rmks}.\!\!}

\Crefname{idea}{{Id}.\!\!}{{Ids}.\!\!}
\newtheorem{warn}[equation]{Warning}
\Crefname{warn}{{Warn}.\!\!}{{Warns}.\!\!}

\Crefname{figure}{{Fig.}\!\!}{{Figs.}\!\!}
\Crefname{footnote}{{Fn.}\!\!}{{Fn.}\!\!}

\Crefname{part}{{\textsection}\!\!}{{\textsection}\!\!}
\Crefname{chapter}{{\textsection}\!\!}{{\textsection}\!\!}
\Crefname{section}{{\textsection}\!\!}{{\textsection}\!\!}
\Crefname{subsection}{{\textsection}\!\!}{{\textsection}\!\!}
\Crefname{appendix}{{\textsection}\!\!}{{\textsection}\!\!}

\usepackage{microtype}

\crefdefaultlabelformat{#2\textup{#1}#3}
\newcommand*{\defeq}{\vcentcolon=}

\newlist{numberenum}{enumerate}{1}
\setlist[numberenum]{label={{\upshape(\arabic*)}}}
\newlist{romanenum}{enumerate}{1}
\setlist[romanenum]{label={\upshape(}\roman*{\upshape)}}
\newlist{letterenum}{enumerate}{1}
\setlist[letterenum]{label={{\upshape(\alph*)}}}

\newcommand*{\pullback}{\arrow[dr, phantom, "\lrcorner", very near start]}
\newcommand*{\Comod}{\mathrm{Comod}}

\newcommand*{\uH}{\mathrm{H}}
\newcommand*{\uE}{\mathrm{E}}

\newcommand*{\D}{\mathcal{D}}
\newcommand*{\calG}{\mathcal{G}}

\newcommand*{\hotimes}{\mathbin{\hat{\otimes}}}

\newcommand{\opgreek}[1]{\begingroup\mathgroup-1 #1\endgroup}

\DeclareMathOperator{\oppi}{\opgreek{\pi}}
\DeclareMathOperator{\opnu}{\opgreek{\nu}}
\DeclareMathOperator{\opsigma}{\opgreek{\sigma}}

\usepackage{titlesec}

\titleformat{\section}
  {\normalfont\Large\bfseries\boldmath}{\thesection}{1em}{}

\titleformat{\subsection}
  {\normalfont\large\bfseries\boldmath}{\thesubsection}{1em}{}

\begin{document}
\title{Descent spectral sequences through synthetic spectra}

\author{Christian Carrick\footnote{\href{mailto:c.d.carrick@uu.nl}{\texttt{c.d.carrick@uu.nl}}},\, Jack Morgan Davies\footnote{\href{mailto:davies@math.uni-bonn.de}{\texttt{davies@math.uni-bonn.de}}},\, and Sven van Nigtevecht\footnote{\href{mailto:s.vannigtevecht@uu.nl}{\texttt{s.vannigtevecht@uu.nl}}}}

\date{July 2, 2025}

\maketitle


\begin{abstract}
The synthetic analogue functor $\nu$ from spectra to synthetic spectra does not preserve all limits.
In this paper, we give a necessary and sufficient criterion for $\nu$ to preserve the global sections of a derived stack. Even when these conditions are not satisfied, our framework still yields synthetic spectra that implement the descent spectral sequence for the structure sheaf, thus placing descent spectral sequences on good footing in the $\infty$-category of synthetic spectra.
As an example, we introduce a new $\mathrm{MU}$-synthetic spectrum $\mathrm{Smf}$.
\end{abstract}

\setcounter{tocdepth}{2}
\tableofcontents


\section*{Introduction}

Real K-theory $\KO$ is the homotopy fixed-points $\KU^{hC_2}$ of complex K-theory under its complex-conjugation action. It is well known that the Adams--Novikov spectral sequence (ANSS) of $\KO$ is isomorphic to the homotopy fixed-point spectral sequence (HFPSS) of $\KU$.
Such identifications are a valuable bridge connecting equivariant and chromatic homotopy theory.
For example, there is an analogous identification of spectral sequences for so-called higher real K-theories $\EO_n$ which plays an essential role in the \emph{Detection Theorem} of Hill--Hopkins--Ravenel \cite{hhr}.\\

The introduction of synthetic spectra by Pstr\k{a}gowski \cite{syntheticspectra} has allowed for a streamlined categorical study of Adams spectral sequences.
In particular, the $\infty$-category $\Syn_\MU$ roughly forms a ``category of Adams--Novikov spectral sequences'' and comes equipped with a functor $\nu\colon\Sp\to\Syn_\MU$, which sends a spectrum to its usual ANSS.
In this article, we use synthetic spectra to prove a much more structured generalization of this identification for a wide range of spectral sequences; see \cref{nuandgacommutesometimes}. 
In the context of the example for $\KO$, this result gives an equivalence of $\E_\infty$-algebras in $\Syn_\MU$ (henceforth \emph{synthetic $\E_\infty$-rings})
\begin{equation}\label{kohfpssandanss}\nu(\KO)\xrightarrow{\cong}\nu(\KU)^{hC_2}.\end{equation}
This is not a tautology, as $\nu$ rarely preserves limits.\\

In this article, we view such an equivalence through the lens of \emph{spectral algebraic geometry}. For example, the \emph{descent spectral sequence} (DSS) for the quotient stack $\Spec (\KU)/C_2$ is equivalent to the HFPSS for $\KU$; we revisit this in \cref{hfpsssignature}. In general, we consider \emph{even-periodic refinements}, which are stacks that are locally constructed from even-periodic affine schemes, such as $\Spec \KU$; see \cite[\textsection2.2]{akhilandlennart} or \cref{evenperiodicrefinenementdefinition} below. In particular, such spectral stacks come with a flat morphism $\x\to \M_\FG$ from the underlying stack to the moduli stack of formal groups.

\begin{theoremalph}[{\cref{theorem:when_Osyn_connective}}]\label{nuandgacommutesometimes}
    Let $(\x,\O^\top_\x)$ be an even-periodic refinement. Suppose that the map $\x \to \M_\FG$ is \emph{tame} {\upshape(}see \cref{tamedefinition}{\upshape)}. Then the natural map of synthetic $\E_\infty$-rings
    \[
        \nu(\O^\top_\x(\x)) \to \lim_{\substack{\Spec A \to \x \\ \text{\upshape\'etale}}} \nu(\O^\top_\x(\Spec A))
    \]
is an equivalence if and only the cohomology groups
\[
    \uH^s(\x \times_{\M_\FG} \Spec L,\, \omega^{\otimes t})
\]
vanish for all $t\in\Z$ and all $s>0$. This condition is satisfied, for example, if $\x\to \M_\FG$ is affine.
\end{theoremalph}

An equivalence as above identifies the ANSS for $\O^\top_\x(\x)$ and the DSS for $\O^\top_\x$ in a highly structured manner. In particular, this recovers unpublished results of Devalapurkar \cite[Th.4]{equivwoods} and Meier \cite[Th.4.6]{relativelyfreetmf} equating these two spectral sequences.\\

When the equivalent conditions of the above theorem fail, the DSS and ANSS for $\O^\top_\x(\x)$ do not usually agree.
In this case, the DSS is often preferable.
We turn this failure of $\nu$ to preserve limits into a definition, resulting in a synthetic spectrum that implements the DSS.
Specifically, we define a \emph{synthetic spectral Deligne--Mumford stack} $(\x,\O^\syn_\x)$ out of the data of an even-periodic refinement $(\x,\O^\top_\x)$; see \cref{definitionofsyntheticlift}.
Its synthetic $\E_\infty$-ring of global sections $\O^\syn_\x(\x)$ is equivalent to the right-hand side of \cref{nuandgacommutesometimes}.
This construction categorifies the DSS, even when the conditions of \cref{nuandgacommutesometimes} fail.

\begin{theoremalph}[{\cref{signatureofsyntheticglobalsections}}]\label{maingeneral}
Let $(\x,\O_\x^\top)$ be an even periodic refinement.
Then the synthetic $\E_\infty$-ring of global sections $\O_\x^\syn(\x)$ implements the descent spectral sequence for the classical global sections $\O^\top_\x(\x)$. Moreover, this identification is one of $\E_\infty$-algebras in filtered spectra.
\end{theoremalph}

An important example is the one of topological modular forms, in both its periodic and projective variants.

\begin{theoremalph}[{\cref{SMFstatements,Smfstatements}}]\label{mainmmf}
The signature spectral sequences of the synthetic $\E_\infty$-rings $\SMF \defeq \O^\syn(\M_\Ell)$ and $\Smf \defeq \O^\syn(\overline{\M}_\Ell)$ are the DSS for $\TMF$ and $\Tmf$, respectively.
There is an equivalence of synthetic $\E_\infty$-rings $\SMF \cong \opnu \TMF$, but the synthetic spectra $\Smf$ and $\opnu \Tmf$ are \textbf{not} equivalent.
\end{theoremalph}

The equivalence $\SMF \cong \opnu\TMF$ follows from \cref{nuandgacommutesometimes}, and recovers the well-known result that the ANSS for $\TMF$ is isomorphic to its DSS.
For $\Tmf$, these spectral sequences differ. In fact, we can quantify the difference between $\Smf$ and $\opnu\Tmf$; see \cref{Smfstatements} for a detailed comparison.\\

The framework of \cref{maingeneral} allows the fruits of synthetic spectra to be brought to bear upon DSSs.
For example, one can now think of the DSS for $\Tmf$ as a modified Adams--Novikov spectral sequence for $\Tmf$.
As a result, tools such as localisations, power operations, and Toda brackets, along with versions of Moss's theorem, may be used more naturally in this context. These ideas are used in \cite{smfcomputation} to calculate the homotopy groups of $\Tmf$ whilst avoiding the circular logic in the current literature. In particular, we highlight that the ANSS for $\Tmf$ is unmanagable and otherwise work with the synthetic spectrum $\Smf$; see \cite[\textsection1.2.1]{smfcomputation}.

\subsection*{Outline}

In \cref{taubsssection}, we collect some facts about $E$-synthetic spectra and their associated spectral sequences, viewed as filtered spectra.
The main result that powers later constructions is \cref{bksssignature}, which gives a way to encode the Bousfield--Kan spectral sequence of a cosimplicial spectrum as an $E$-synthetic spectrum.
We also review the spectral sequence associated to the synthetic analogue of a spectrum, identifying it with the $E$-Adams spectral sequence; see \cref{signatureofgammanu}.
Except for \cref{bksssignature}, this material is well known; we review it to establish notation to be used later on.\\


From \cref{sec:Osyn} onwards, we specialise to $\MU$-synthetic spectra.
In \cref{sec:Osyn}, we show that the synthetic global sections $\O^\syn_\x(\x)$ of an even-periodic stack $(\x,\O^\top_\x)$ implement the DSS for $\O^\top_\x$, proving \cref{maingeneral}. We show that the synthetic spectrum $\O^\syn_\x(\x)$ has connective cover $\opnu(\O^\top_\x(\x))$ (\cref{invert_tau_Osyn}), so that studying the difference comparison of the DSS for $\O^\top_\x$ and the ANSS for $\O^\top_\x(\x)$ comes down to studying synthetic connectivity of $\O^\syn_\x(\x)$.\\

In \cref{sec:computing_connectivity}, we compute the connectivity of the global sections $\O^\syn_\x(\x)$.
The connectivity of a general $\MU$-synthetic spectrum $S$ is controlled by the homology $\nu\MU_{*,*}(S)$, and accordingly, the goal of this section is to set up a spectral sequence to compute this homology.
The difficult part is identifying the abutment of this spectral sequence; we do this in the more technical \cref{ssec:identifying_global_sections}.
Using this, in \cref{sec:comparison_DSS_ANSS} we prove \cref{nuandgacommutesometimes}, or a more precise version in \cref{theorem:when_Osyn_connective}.\\

In \cref{sec:application}, we apply the above theory in the two examples mentioned in the introduction.
First, in \cref{sec:HFPSS}, we consider the fixed-points of the action of a finite group $G$ on the a Landweber exact $\E_\infty$-ring.
This includes the action of a finite subgroup of the Morava stabiliser group on Morava E-theory, and also the case of K-theory mentioned in the introduction.
Second, we study topological modular forms in \cref{sec:SMF}.
We introduce the synthetic spectra $\SMF$ and $\Smf$ in \cref{smfdefinition}, and then show that $\opnu\TMF\cong \SMF$ and $\opnu\Tmf\ncong \Smf$, proving \cref{mainmmf}.

\subsection*{Conventions}
\begin{itemize}
    \item The language of $\infty$-categories \`{a} la \cite{httname} and \cite{sagname} will be used. Our limits and colimits will only be taken over small $\infty$-categories with respect to some fixed Grothendieck universe.
    \item What Lurie calls a nonconnective qc separated spectral Deligne--Mumford stack whose underlying $\infty$-topos is $1$-localic, we will simply refer to as a \emph{spectral Deligne--Mumford stack}.
    \item We assume that all classical Deligne--Mumford stacks are qc and separated.
    Consequently, by \cite[Lm.B.2]{tmfwls}, we can write a spectral Deligne--Mumford stack as a pair $(\x,\O^\top_\x)$, where $\x$ is a classical Deligne--Mumford stack, and where $\O^\top_\x$ is an \'{e}tale hypersheaf of $\E_\infty$-rings on the small \'{e}tale site of $\x$.
    \item The homotopy groups of a sheaf of $\E_\infty$-rings will always implicitly be sheafified.
    \item The small \'{e}tale site of a Deligne--Mumford stack $\x$ will be denoted by $\DM^\et_{/\x}$, and we write $\Aff^\et_{/\x}$ for the subsite spanned by affine stacks.
    \item Spectral sequences will always use Adams indexing, meaning that in $\mathrm{E}_r^{k,s}$, the variable $k$ denotes the stem, $s$ denotes the filtration, and the differential $d_r$ has bidegree $(-1,r)$.
\end{itemize}

\subsection*{Acknowledgements}

We thank Shaul Barkan, Sil Linskens, Lennart Meier, and Lucas Piessevaux for their stimulating discussions around these topics, and Lennart Meier for his feedback on earlier drafts of this article. Thank you as well to the anonymous referee for their suggestions.\\

The first author was supported by NSF grant \texttt{DMS-2401918}.
The first and third author were supported by the NWO grant \texttt{VI.Vidi.193.111}.
The second author is an associate member of the Hausdorff Center for Mathematics at the University of Bonn (\texttt{DFG GZ 2047/1}, project ID \texttt{390685813}).


\section{Synthetic spectra and spectral sequences}\label{taubsssection}

For this section, let $E$ be a homotopy associative ring spectrum of Adams type (\cite[Df.3.14]{syntheticspectra}), and write $\Syn_E$ for the $\infty$-category of $E$-based synthetic spectra of \cite{syntheticspectra}. For an introduction to synthetic spectra, we refer to \cite{syntheticspectra} and \cite[\textsection9]{burkhahnseng}, and we further follow much of the notation from \cite[\textsection2--3]{syntheticj}.

\begin{notation}
    We use ``stem--filtration'' grading for synthetic spectra, so that $\pi_{k,s}$ corresponds to a $(k,s)$-location in an Adams chart.
    Formally, this means that we write $\Sph^{a,b}\defeq\Sigma^{-b}\nu\Sph^{a+b}$.
    In particular, $\tau$ has bidegree $(0,-1)$ and the $\infty$-categorical suspension has bidegree $(1,-1)$.
    Note that this differs from the grading used in \cite{syntheticspectra}.
\end{notation}

We will refer to an $\E_\infty$-object in $\Syn_E$ as a \emph{synthetic $\E_\infty$-ring}.

\begin{notation}
    The stable $\infty$-category $\Syn_E$ has a natural t-structure \cite[Pr.4.16]{syntheticspectra}.
    Following \cite{burkhahnseng}, we call this the \emph{homological t-structure}. In this t-structure, a synthetic spectrum $X$ is $n$-connective if and only if $\nu E_{\ast,s}(X)=0$ for all $s>-n$; see \cite[Th.4.18]{syntheticspectra}.
\end{notation}


\subsection{The synthetic analogue functor}
\label{ssec:nu}

In this subsection we list a few standard facts about the \emph{synthetic analogue} functor $\nu\colon \Sp\to \Syn_E$ of \cite[Df.4.3]{syntheticspectra}.

\begin{prop}
    \label{prop:nu_of_E_module}
    \leavevmode
    \begin{numberenum}
        \item Let $M$ be a spectrum that admits the structure of a homotopy $E$-module.
        Then we have a natural isomorphism of bigraded abelian groups
        \[
            \pi_{*,*}\nu M \cong \pi_* M [\tau],
        \]
        where $\pi_k M$ is placed in bidegree $(k,0)$.
        \item \label{item:nu_strong_monoidal} Let $X$ be spectrum that can be written as a filtered colimit of finite $E$-projective spectra.
        Then for any spectrum $Y$, the natural map
        \[
            \nu X \otimes \nu Y \to \nu(X \otimes Y)
        \]
        coming from the lax monoidal structure on $\nu$ is an equivalence.
        In particular, this holds if $X = E$.
    \end{numberenum}
\end{prop}
\begin{proof}
    The first item is \cite[Pr.4.60]{syntheticspectra}, and the second is \cite[Lm.4.24]{syntheticspectra}.
\end{proof}

\begin{example}\label{comeansmumodule}
    As in most of this paper, consider the case where $E = \MU$, so that $\Syn_\MU$ is the categorification of the Adams--Novikov spectral sequence.
    If $R$ is a Landweber-exact ring spectrum, then $R$ can be written as a filtered colimit of finite even spectra; see \cite[§2.1]{hoveysticklandMoravalocal}.
    A finite even spectrum in particular has projective $\MU_*$-homology, so \cref{prop:nu_of_E_module}\,\ref{item:nu_strong_monoidal} applies.
\end{example}

Next, we discuss the relationship between $\nu$ and the $\tau$-inversion functor $\tau^{-1} \colon \Syn_E \to \Sp$ of \cite[Df.4.39]{syntheticspectra}.
Let us introduce some terminology.

\begin{mydef}
    Let $X$ be a spectrum.
    A \emph{synthetic lift} of $X$ is a synthetic spectrum $S$ such that $\tau^{-1} S$ is equivalent to $X$.
\end{mydef}

The functor $\tau^{-1}$ is a left inverse to $\nu$ by \cite[Pr.4.40]{syntheticspectra}, so the synthetic analogue is an example of a synthetic lift. More generally, the following fact often allows for the construction of a much broader class of lifts.

\begin{prop}
    \label{prop:limits_are_synthetic_lifts}
    Let $X \colon I \to \Sp$ be a diagram of spectra.
    Then we have an equivalence
    \[
        \nu\left(\lim_{\alpha \in I}X_\alpha\right) \cong \tau_{\geq 0} \left(\lim_{\alpha \in I} \nu (X_\alpha)\right),
    \]
    and the natural limit-comparison map
    \[
        \tau^{-1}\left(\lim_{\alpha \in I}\nu(X_\alpha)\right) \to \lim_{\alpha \in I} X_\alpha
    \]
    is an equivalence.
    In particular, the synthetic spectrum $\lim_{\alpha \in I}\nu(X_\alpha)$ is a synthetic lift of the spectrum $\lim_{\alpha\in I}X_\alpha$.
\end{prop}

\begin{proof}
    Recall from \cite[Pr.4.21]{syntheticspectra} that $\nu$ lands in connective synthetic spectra, so that we can restrict it to a functor $\Sp \to \Syn_E^{\geq 0}$. By \cite[Pr.4.36]{syntheticspectra}, the functor $\nu$ is then equivalent to the composite
    \[
        \Sp \cong \Syn_E[\tau^{-1}] \hookrightarrow \Syn_E \to \Syn_E^{\geq 0}
    \]
    of the equivalence between spectra and $\tau$-invertible synthetic spectra \cite[Th.4.37]{syntheticspectra}, the inclusion into all synthetic spectra, and the connective cover functor.
    By \cite[Pr.4.33]{syntheticspectra}, the inclusion $\Syn_E[\tau^{-1}] \hookrightarrow \Syn_E$ is right adjoint to $\tau^{-1} \colon \Syn_E\to \Syn_E[\tau^{-1}]$.
    We see that $\nu \colon \Sp \to\Syn_E^{\geq 0}$ is right adjoint to the composite
    \[
        \Syn_E^{\geq 0} \hookrightarrow \Syn_E \to \Syn_E[\tau^{-1}] \cong \Sp
    \]
    of the inclusion and $\tau$-inversion, i.e., $\nu \colon \Sp \to \Syn_E^{\geq0}$ is right adjoint to $\tau^{-1} \colon \Syn_E^{\geq0} \to \Sp$.\\
    
    In particular, $\nu$ sends a limit $\lim_{\alpha\in I} X_\alpha$ of spectra to the synthetic spectrum
    \[
        \tau_{\geq 0}\left(\lim_{\alpha\in I} \nu(X_\alpha)\right).
    \]
    Since $\tau^{-1}$ is a left inverse to $\nu$, it suffices now to observe that if $S$ is any synthetic spectrum, then the connective cover $\tau_{\geq 0}S \to S$ becomes an equivalence upon inverting $\tau$.
    This fact follows from \cite[Lm.4.35]{syntheticspectra}, which says that any bounded-above synthetic spectrum becomes zero after inverting~$\tau$.
\end{proof}

\begin{remark}
    Because $\nu$ is fully faithful, the above corollary in fact says that $\tau^{-1}$ preserves the limit of any diagram that takes values in synthetic analogues.
\end{remark}

\begin{warn}
    The functor $\tau^{-1}$ does not necessarily preserve limits when the values of the diagram are not synthetic analogues.
    For example, if $X$ is an $E$-nilpotent complete spectrum, then $\nu X$ is the limit of its $\tau$-adic tower $\{ \nu X/\tau^n \}_n$, but every term in this tower becomes zero after inverting $\tau$.
\end{warn}


\subsection{Filtered spectra}
\label{sec:syn_and_sseq}
In this paper, we view the $\infty$-category $\Fil(\Sp)=\Fun(\Z^\op,\Sp)$ of \emph{filtered spectra} as a more structured version of a `category of spectral sequences'.
For an explanation of how a filtered spectrum gives rise to a spectral sequence, see, e.g., \cite[\textsection 1.2.2]{haname} or \cite[\textsection II]{alicethesis}, or \cite[\textsection2]{syntheticj} for a summary more aligned with the notation we use here.
The additional structure of a filtered spectrum allows one, for instance, to talk about coherently commutative multiplicative structures on a spectral sequence.\\

Meanwhile, for $E$ an Adams-type homology theory, the $\infty$-category $\Syn_E$ should be thought of as an $\infty$-category of $E$-based Adams spectral sequences.
These perspectives fit together via a functor $\Syn_E \to \Fil(\Sp)$ which we now recall.
This construction is well-known: in \cite[Def.5.58]{iraklipiotr} this is used on synthetic analogues to construct the Adams filtration (see \cref{ssec:signature_nu} below for a further discussion).
The following form of this construction we recall from \cite{syntheticj}.

\begin{mydef}
    Let $\sigma\colon \Syn_E\to \Fil(\Sp)$ be the functor given by applying the mapping spectrum functor $\map(\Sph,{-})$ level-wise to the $\tau$-tower of a synthetic spectrum.
    More explicitly, the map $\sigma(X)^n \to \sigma(X)^{n-1}$ is of the form
    \[\sigma(X)^n=\map_{\Syn_E}(\Sph,\, \Sigma^{0,-n}X)\xrightarrow{\map(\Sph,\tau)}\map_{\Syn_E}(\Sph,\, \Sigma^{0,-n+1}X)=\sigma(X)^{n-1}.\]
    The functor $\sigma$ comes equipped with a natural lax symmetric monoidal structure: see \cite[Th.2.5(4)]{syntheticj}, where the notation $\Gamma$ is used for $\sigma$.
\end{mydef}

\begin{mydef}\label{signaturedefn}
    Let $X$ be a synthetic spectrum. Then the \emph{signature} of $X$ is the filtered spectrum $\sigma X$ and the associated spectral sequence is called the \emph{signature spectral sequence}.
\end{mydef}

The functor $\sigma$ is a right adjoint\footnote{The left adjoint to $\sigma$ is the universal extension of the functor $\Z \to \Syn_E$ defined by the $\tau$-tower of the unit.} and, in many cases, part of a monadic adjunction \cite[Cor.6.1]{lawsoncellular}. For now, we will just need the following.

\begin{lemma}\label{gammaandlimits}
    The functor $\sigma$ preserves limits.
\end{lemma}

\begin{proof}
    This is clear from the universal property of limits and the formula for $\sigma$:
    \[\sigma(\lim X_\al)^n=\map_{\Syn_E}(\Sph,\, \Sigma^{0,-n}\lim X_\al)\cong \lim \map_{\Syn_E}(\Sph,\, \Sigma^{0,-n} X_\al)=\lim \sigma(X_\al)^n. \qedhere\]
\end{proof}

In \cite{syntheticj}, the first two authors study spectral sequences arising in this way which have no obvious other interpretation.
In the remainder of this section on the other hand, we want to identify these spectral sequences with a variety of classical ones.


\subsection{Bousfield--Kan spectral sequences}

The goal of this subsection is to construct a filtered spectrum that encodes the Bousfield--Kan spectral sequence (BKSS) of a cosimplicial spectrum.
In practise, such a cosimplicial spectrum comes about as a choice of resolution of a certain spectrum.
The $\uE_1$-page of the BKSS depends on this choice, but the $\uE_2$-page does not.
To ensure that our filtered spectra do not depend on this choice of resolution, we will apply a \emph{d\'ecalage} functor by definition.
This construction is well-known, and is essentially the same as the one in \cite[\textsection3]{mmf}.
We end this subsection with a new result on when these filtered spectra can be lifted to an $E$-based synthetic spectrum; see \cref{bksssignature}.\\

To define the d\'ecalage, let us recall a definition of the Whitehead tower functor $\Wh\colon \Sp\to \Fil(\Sp)$.
There is a t-structure on $\Fil(\Sp)$ whose connective objects are those filtered spectra $X$ such that $X^n$ is $n$-connective for each $n\in \Z$; this is called the \emph{canonical t-structure} by Hedenlund \cite[Pr.II.1.22]{alicethesis}. Moreover, this t-structure is compatible with the symmetric monoidal structure on $\Fil(\Sp)$ given by Day convolution; see \cite[Pr.II.1.23]{alicethesis}.
As a result, the connective cover functor $\tau_{\geq 0}$ is lax symmetric monoidal.

\begin{mydef}\label{defofwhitehead}
    Define a functor ${\WH}\colon \Sp\to \Fil(\Sp)$, the \emph{Whitehead tower} functor, to send a spectrum $X$ to the filtered spectrum $\tau_{\geq 0}\Const( X)$ where $\Const\colon \Sp\to \Fil(\Sp)$ is the constant functor. More explicitly, $\Wh(X)$ can be written as
    \[\cdots \to \tau_{\geq n+1}X \to \tau_{\geq n}X\to \tau_{\geq n-1}X\to \cdots.\]
    By \cite[Pr.II.1.26]{alicethesis}, and the arguments above, the functor $\Wh$ is lax monoidal.
\end{mydef}

Using this, we can define the filtered spectrum encoding the BKSS.
In what follows, and the remainder of this article, we give the $\infty$-category $\Sp^\Delta$ of cosimplicial spectra the levelwise symmetric monoidal structure.

\begin{mydef}\label{definitionofdecalage}
    Define the \emph{cosimplicial d\'{e}calage} functor ${\Dec}\colon \Sp^\Delta \to \Fil(\Sp)$ by the formula
    \[\Dec(X^\bullet)=\lim_\Delta (\WH(X^\bullet)).\]
    Notice $\Dec$ is lax symmetric monoidal as it is the composition
    \[\Sp^\Delta = \Fun(\Delta, \Sp) \xrightarrow{\Fun(\Delta, \Wh)} \Fun(\Delta, \Fil(\Sp)) \xrightarrow{\lim_\Delta} \Fil(\Sp)\]
    where \cref{defofwhitehead} states the first functor is lax symmetric monoidal, and the second is too as a limit.
\end{mydef}

\begin{remark}
    \label{rmk:Dec_semicosimplicial}
    Since totalisation only depends on the underlying semicosimplicial object, the functor $\Dec$ factors through the forgetful functor $\Sp^\Delta \to \Sp^{\Delta_{\mathrm{inj}}}$ to semicosimplicial spectra.
\end{remark}

In order to ensure (conditional) convergence, it is important to note that these filtered spectra are complete in the following sense.

\begin{prop}
    \label{lem:decalage_is_complete}
    Let $X^\bullet$ be a cosimplicial spectrum.
    Then the filtered spectrum $\Dec(X^\bullet)$ is a \emph{complete filtration} of $\lim_\Delta X^\bullet$, meaning that $\lim_k \Dec(X^\bullet)^k$ vanishes and $\colim_k \Dec(X^\bullet)^k \cong \lim_\Delta X^\bullet$.
\end{prop}

\begin{proof}
    The limit of the Whitehead tower of a spectrum vanishes, so the claim about the limit follows since limits commute with limits.
    
    To see that $\colim \Dec(X^\bullet)\cong \lim_\Delta X^\bullet$,  note that, for any $k\in\Z$, one has a fibre sequence
    \[\tau_{\ge k+1}X^\bullet\to X^\bullet\to \tau_{\le k}X^\bullet.\]
    Applying $\lim_{\Delta}$, one therefore 
    has a fibre sequence
    \[\lim_{\Delta}\tau_{\ge k+1}X^\bullet\to \lim_{\Delta}X^\bullet\to \lim_{\Delta}\tau_{\le k}X^\bullet.\]
    Applying $\colim_k$, one has a cofibre sequence
    \[\colim_k\lim_{\Delta}\tau_{\ge k+1}X^\bullet\to \lim_{\Delta}X^\bullet\to \colim_k\lim_{\Delta}\tau_{\le k}X^\bullet.\]
    By definition, the lefthand term gives the colimit of $\Dec(X^\bullet)$, so it suffices to show the righthand term vanishes. Since coconnectivity is preserved by limits, we see that for all $k\in\Z$, the spectrum $\lim_{\Delta}\tau_{\le k}X^\bullet$ is $k$-truncated. Letting $k$ tend to $\infty$, we obtain a spectrum that is $k$-truncated for all $k$, which is therefore zero.
\end{proof}

The relationship between this definition of d\'{e}calage and the ``page-turning'' intuition is given by the following.

\begin{prop}\label{decalageidentification}
    Let $X^\bullet$ be a cosimplicial spectrum, and let $\{\mathrm{E}_r^{k,s}(X^\bullet)\}_{r\geq 1}$ denote the Bousfield--Kan spectral sequence for $X^\bullet$.
    Then there is an isomorphism of spectral sequences {\upshape(}where $r\geq 2${\upshape)}
    \[
        \mathrm{E}_r^{k,s}(\Dec(X^\bullet)) \cong \mathrm{E}_{r+1}^{k,\, s-k}(X^\bullet).
    \]
\end{prop}
\begin{proof}
    See \cite[Pr.6.3]{levinedecalage}, but note that Levine uses a different indexing from the Adams indexing used above.
\end{proof}

\begin{notation}
    The spectral sequence associated to a filtered spectrum is usually indexed to start on the $\uE_1$-page.
    Using the reindexing of the previous result, we will instead reindex the spectral sequence arising from the décalage construction of \cref{definitionofdecalage} to begin on the $\uE_2$-page.
\end{notation}

In particular, if $X^\bullet$ is a cosimplicial spectrum, then the BKSS computing the homotopy groups of the limit of $X^\bullet$ is the spectral sequence associated with the filtered spectrum $\Dec(X^\bullet)$. If the resolution $X^\bullet$ is nice enough, then we can package this BKSS as a synthetic spectrum, or more precisely as a synthetic lift of $\lim X^\bullet$.
Moreover, this synthetic spectrum recovers the filtered spectrum $\Dec(X^\bullet)$, together with its multiplicative structure. A similar results appears in \cite[Lm.9.18]{barthelpstragowski}.

\begin{theorem}\label{bksssignature}
    Let $X^\bullet$ be a cosimplicial spectrum, and write $S$ for the limit of the cosimplicial synthetic spectrum $\nu X^\bullet$. If each $X^i$ admits the structure of a homotopy $E$-module, then there is an equivalence of filtered spectra
    \[\sigma(S)\cong \Dec(X^\bullet).\]
    In particular, if each $X^i$ admits a homotopy $E$-module structure, then $S$ has the signature of the BKSS for $X^\bullet$.
    Moreover, the above equivalence is natural in the following way: it is a symmetric monoidal equivalence between the lax symmetric monoidal functors $\sigma \circ {\lim_\Delta} \circ \nu$ and $\Dec$ when restricted to the cosimplicial spectra of the above kind.
\end{theorem}

For the proof, we require the following lemma.
Note that a similar technique is employed to prove \cite[Pr.C.22]{burkhahnsenggalois}.

\begin{lemma}\label{homotopymodulewhitehead}
    Fix a homotopy associative ring spectrum $E$ of Adams type.
    Let $X$ be a spectrum admitting a homotopy $E$-module structure.
    Then there is an equivalence of filtered spectra
    \[\sigma(\nu X)\cong \WH(X).\]
    This equivalence is natural in the following sense: it is a symmetric monoidal equivalence between the lax symmetric monoidal functors $\sigma\circ \nu$ and $\Wh$ when restricted to spectra admitting a homotopy $E$-module structure.
\end{lemma}

\begin{proof}
    For each $n\in \Z$ there is a natural map of spectra
    \[
        (\sigma \nu X)^n=\map_{\Syn_E}(\Sph,\, \Sigma^{0,-n}\nu X)\xrightarrow{\tau^{-1}} \map_{\Sp}(\Sph,X)\cong X
    \]
    induced by $\tau$-inversion.
    This lifts to a monoidal natural transformation $\sigma(\nu X)\to \Const(X)$ as $\tau$-inversion itself is symmetric monoidal; see \cite[Pr.4.40]{syntheticspectra}.
    We now observe that this transformation naturally factors through the connective cover of $\Const(X)$ with respect to the canonical t-structure discussed prior to \cref{defofwhitehead}.
    Indeed, to calculate
    \[\pi_k (\sigma \nu X)^n \cong \oppi_0 \map_{\Syn_E}(\Sph^{k,\, n-k},\,\nu X)=\oppi_{k,\, n-k}\nu X\]
    we appeal to \cite[Pr.4.60]{syntheticspectra}, which states that the group $\pi_{k,s}\nu X$ vanishes for $s>0$, so the above group vanishes for $k<n$. This factorisation yields our desired natural symmetric monoidal comparison map
    \[\sigma(\nu X) \to \tau_{\geq 0} (\mathrm{Const}(X)) = \Wh(X).\]
    To see this map is an equivalence, it suffices to show that $\sigma(\nu X)^n \to \tau_{\geq n}X$ is an equivalence for all $n$, and this follows from another application of \cite[Pr.4.60]{syntheticspectra}.
\end{proof}

\begin{proof}[Proof of \cref{bksssignature}]
    We have the natural monoidal equivalences of filtered spectra
    \[\Dec(X^\bullet)=\lim( \WH(X^\bullet)) \cong \lim \sigma(\nu X^\bullet)\cong \sigma \lim (\nu X^\bullet) = \sigma(S)\]
    from \cref{homotopymodulewhitehead,gammaandlimits}, respectively.
\end{proof}

Observe that \cref{prop:limits_are_synthetic_lifts} implies that the synthetic spectrum $S$ in the above theorem is a synthetic lift of $\lim X^\bullet$.

\begin{remark}
    If $X^\bullet$ is an $\E_n$-algebra in $\Sp^\Delta$, then the above proposition in particular results in an $\E_n$-structure on the synthetic spectrum $S$ as in the proposition.
    This stands in contrast to an important application of synthetic spectra, in that they allow for multiplicative structures even when cosimplicial models are not known to have such structures; see the discussion in \cref{ssec:signature_nu}.
    Nevertheless, one of the main perspectives of this article is that a synthetic refinement of a filtered spectrum has many benefits, such as allowing for a comparison of its underlying spectral sequence to Adams spectral sequences.
\end{remark}

\subsection{The signature of a synthetic analogue}
\label{ssec:signature_nu}

Classically, one can define a filtered model of the Adams filtration for any spectrum $E$ with a map $\Sph \to E$; see \cref{rmk:filtered_models_using_less_structure} below.
However, it is difficult to upgrade this functor to a lax monoidal functor.
As recalled in, e.g., \cite[\textsection5.5]{iraklipiotr}, if $E$ has the structure of an $\E_\infty$-ring, then this classical definition of the Adams filtration can also be made into a lax symmetric monoidal functor.
The $E$-Adams spectral sequence itself however does not depend on this multiplicative structure on $E$.
This makes it infeasible as a general approach to define a monoidal Adams filtration functor, as one would also like to do this for $E$ that cannot be given an $\E_\infty$-structure.

One of the great breakthroughs of synthetic spectra is that they allow one to define a lax monoidal structure on the Adams filtration in the most general case.
Specifically, Patchkoria and Pstr\k{a}gowski \cite[\textsection5.4]{iraklipiotr} show that (to use the language of the present paper) the filtered spectrum $\sigma\nu(X)$ captures the Adams filtration on $\pi_*(X)$; see \cite[Th.5.60]{iraklipiotr}.
As the functor $\sigma \circ \nu$ is canonically lax symmetric monoidal, this results in a coherently monoidal version of the Adams filtration.
Importantly, this construction requires barely any multiplicative structure on~$E$: the $\infty$-category $\Syn_E$ is symmetric monoidal for all homotopy associative ring spectra $E$ of Adams type.

\begin{remark}
More precisely, Patchkoria and Pstr\k{a}gowski \cite{iraklipiotr} introduce and use a more general type of synthetic category than the one originally introduced in \cite{syntheticspectra}.
This is done to allow for much more general $E$ than homotopy ring spectra of Adams type.
In the end, they construct a lax symmetric monoidal functor $\Sp \to \Fil(\Sp)$ encoding the $E$-Adams filtration whenever $E$ is a spectrum that admits right-unital multiplication; see \cite[Th.5.64]{iraklipiotr}.
\end{remark}

For most of this article, we find ourselves in the situation where $E$ has an $\E_\infty$-ring structure: our comparison results happen in the case where $E = \MU$.
For the sake of completeness, we will in this section compare this classical definition of the Adams filtration with the modern definition $\sigma \circ \nu$, and show that the monoidal structures also agree.
This is certainly well known: see \cite[Rmks.4.62--65]{syntheticspectra}, \cite[\textsection A.1]{burkhahnseng}, and \cite[Df.5.58, Th.5.60]{iraklipiotr}, but we include it for completeness and to establish notation to be used later.

\begin{mydef}\label{adamstowerdefinition}
    Let $E$ be an $\E_1$-ring.
    Via the construction of \cite[Constr.2.7]{mnn}, the $\E_1$-structure on $E$ gives rise to a cosimplicial spectrum $E^{\otimes \bullet+1}$.
    Define the functor $\ASS_E\colon \Sp\to \Sp^{\Delta}$ to send a spectrum $X$ to its \emph{cosimplicial Adams tower}
    \[\ASS_E(X)^\bullet \defeq X\otimes E^{\otimes \bullet+1}.\]
    If $E$ moreover has an $\E_\infty$-structure, then by \cite[Lm.2.39]{achimthesis} this functor has a canonical lax symmetric monoidal structure.
\end{mydef}

\begin{remark}
    If $E$ is an $\E_k$-ring for $k\geq 1$, then the functor $\ASS_E$ is lax $\E_k$-monoidal by the same \cite[Lm.2.39]{achimthesis}.
    In what follows, we will for simplicity focus only on the case where $E$ is an $\E_\infty$-ring, but all results hold when $E$ is an $\E_k$-rings as well, provided ``symmetric monoidal'' is replaced by ``$\E_k$-monoidal''.
\end{remark}

\begin{remark}
    \label{rmk:filtered_models_using_less_structure}
    As noted before, one does not need an $\E_1$-structure on $E$ to define a filtered model for the $E$-Adams spectral sequence.
    Indeed, to construct the \emph{semicosimplicial} object underlying $X\otimes E^{\otimes \bullet+1}$, one only needs the map $\Sph \to E$, without any other multiplicative structure on $E$.
    The functor $\Dec$ only depends on the underlying \emph{semicosimiplicial} object (see \cref{rmk:Dec_semicosimplicial}), so applying it to $X\otimes E^{\otimes\bullet+1}$ results in a functorial filtered model for the Adams spectral sequence.
    Alternatively, as noted in \cite[Rmk.2.27]{iraklipiotr}, one can also define a \emph{cosimplicial} spectrum giving rise to the $E$-Adams spectral sequence for $X$, without any multiplicative structure on $E$ (as explained in loc.\ cit., this can be done for any Adams spectral sequence in the sense of Miller).
    However, such a cosimplicial spectrum will be far from canonical, and importantly, will not be functorial in $X$; cf.\ \cite[Rmk.5.62]{iraklipiotr}.
    An $\E_1$-structure on $E$ meanwhile allows one to construct the particular cosimplicial tower of \cref{adamstowerdefinition}, which is functorial in $X$ by construction.
\end{remark}

If $E$ is an $\E_\infty$-ring, then the functor $\Dec \circ \ASS_E$ is a lax symmetric monoidal functor.
We can then compare the lax symmetric monoidal functors $\Dec \circ \ASS_E$ and $\sigma \circ \nu$.

\begin{prop}\label{signatureofgammanu}
    Fix a homotopy associative ring spectrum $E$ of Adams type.
    Let $X$ be an $E$-nilpotent complete spectrum.
    Then there is an equivalence of filtered spectra
    \[\sigma\nu (X)\cong \Dec(\ASS_E(X)).\]
    This equivalence is natural in the following sense: it is an equivalence between the functors $\sigma\circ \nu$ and $\Dec\circ \ASS_E$ when restricted to $E$-nilpotent complete spectra.
    If $E$ is moreover an $\E_\infty$-ring, then this equivalence is canonically symmetric monoidal.
\end{prop}

\begin{proof}
    Consider the natural (in $X$) equivalences of filtered spectra
    \begin{equation}\label{gammatotofantower}
        \Dec(\ASS_E(X))=\lim \left(\WH(X\otimes E^{\otimes \bullet+1})\right) \cong \lim \sigma \left(\nu (X\otimes E^{\otimes\bullet+1})\right) \cong  \sigma \lim\left(\nu (X\otimes E^{\otimes \bullet+1})\right)
    \end{equation}
    from \cref{homotopymodulewhitehead} and \cref{gammaandlimits}, respectively.
    By \cref{prop:nu_of_E_module}, for any $k\geq 1$ we can identify $\nu(X\otimes E^{\otimes k})$ with $\nu X \otimes (\nu E)^{\otimes k}$, and hence identify \eqref{gammatotofantower} with $\sigma\lim ((\nu X)\otimes (\nu E)^{\otimes \bullet})$.
    As the $\nu E$-nilpotent completion of $\nu X$ is $\nu$ applied to the $E$-nilpotent completion of $X$ by \cite[Pr.A.13]{burkhahnseng}, we see that \eqref{gammatotofantower} is naturally equivalent to
    \[\sigma \lim\left(\nu (X\otimes E^{\otimes \bullet+1})\right)\cong \sigma\lim \left(\nu X\otimes (\nu E)^{\otimes \bullet+1}\right)\cong \sigma \nu(X^{\wedge}_E)\cong \sigma (\nu X).\]
    If $E$ is an $\E_\infty$-ring, then the equivalences of \eqref{gammatotofantower} are symmetric monoidal in $X$.
\end{proof}


\section{The synthetic lift of an even-periodic refinement}
\label{sec:Osyn}

If $\x$ is a stack equipped with a sheaf of spectra $\F$, then the \emph{descent spectral sequence} (DSS) for $\F$ attempts to compute the homotopy groups of the global sections $\Gamma(X,\F)$ by means of the homotopy sheaves $\pi_* \F$. We study the DSS from a synthetic point of view in the context of an \emph{even-periodic refinement} of a classical Deligne--Mumford stack.
We show that, in this case, the DSS is most naturally viewed as the signature spectral sequence of a certain $\MU$-synthetic spectrum.
First, we define a filtered model for the DSS in \cref{sec:DSS}, and then construct this $\MU$-synthetic spectrum in \cref{ssec:even_periodic_refinement}.

\subsection{Descent spectral sequences in filtered spectra}\label{sec:DSS}
We will consider the descent spectral sequence for sheaves $\F$ of spectra defined on the small \'{e}tale site of a classical Deligne--Mumford stack $\x$.

\begin{construction}\label{dssdefinition}
    Let $\x$ be a classical Deligne--Mumford stack, $\Spec A\to \x$ an affine \'{e}tale cover of $\x$, and $\F$ be a sheaf\footnote{All of the results of this paper hold if one only assumes our spectral Deligne--Mumford stacks are \emph{quasi-separated} rather than separated, if all sheaves are taken to be hypersheaves.} of spectra on the small \'{e}tale site of $\x$.
    Let $\Spec A^\bullet\to \x$ be the \v{C}ech nerve associated to the cover $\Spec A \to \x$;\footnote{Note that such an affine \'{e}tale cover and an affine \'{e}tale \v{C}ech nerve exist as we assume our Deligne--Mumford stacks are qc and separated, respectively.} note that this is an augmented simplicial object over $\x$.
    We define the \emph{descent spectral sequence} (DSS) for the pair $(\x,\F)$ as the filtered spectrum
    \[\DSS(\x, \F)=\Dec(\F(\Spec A^\bullet)).\]
    As $\Dec$ is lax symmetric monoidal by \cref{definitionofdecalage}, the filtered spectrum $\DSS(\x,\F)$ inherits any $\E_n$-structure from $\F$. As $\F$ is an \'{e}tale sheaf, the spectral sequence associated to this filtered spectrum then abuts to the global sections of $\F$:
    \[\Ga(\x,\F)\xrightarrow{\cong} \lim_\Delta \F(\Spec A^\bullet).\]
\end{construction}

The following two lemmata justify this definition.

\begin{lemma}\label{e2pageofdssidentification}
    Let $(\x,\F)$ be as in \cref{dssdefinition}. Then the $\mathrm{E}_2$-page of the spectral sequence associated to $\DSS(\x,\F)$ takes the form
    \[\mathrm{E}_2^{k,s} \cong \mathrm{H}^s(\x,\pi_{k+s}\F).\]
\end{lemma}

\begin{proof}
    The complex computing the sheaf cohomology of $\pi_{k+s} \F$ can be constructed taking an affine \'{e}tale cover of $\x$ and then evaluating $\F$ on the \v{C}ech nerve of this cover. By definition, this is precisely the $\mathrm{E}_1$-page of the spectral sequence associated to $\DSS(\x,\F)$, and the result follows.
\end{proof}

\begin{lemma}
    \label{lem:DSS_independent_of_resolution}
    Let $(\x,\F)$ be as in \cref{dssdefinition}, choose two affine \'{e}tale covers $\Spec A_i\to \x$ of $\x$ for $i=1,2$, and write $\DSS_i(\x,\F)$ for the associated filtered spectra \`{a} la \cref{dssdefinition}. Then there is a preferred monoidal equivalence of filtered spectra
    \[\DSS_1(\x,\F) \cong \DSS_2(\x,\F),\]
    meaning a chosen equivalence, which is one of $\E_n$-objects if $\F$ has an $\E_n$-structure.
\end{lemma}

\begin{proof}
    Given two affine \'{e}tale covers of $\x$ with \v{C}ech nerves $\Spec A_1^\bullet$ and $\Spec A_2^\bullet$, we define a third affine \'{e}tale cover $\Spec B$ as the common pullback of stacks:
    \[\begin{tikzcd}
        {\Spec B}\ar[d]\ar[r] \arrow[dr, phantom, "\lrcorner", very near start]  &   {\Spec A_1}\ar[d]   \\
        {\Spec A_2}\ar[r]   &   {\x}.
    \end{tikzcd}\]
    We also have an identification of \v{C}ech nerves
    \[\Spec B^\bullet \cong \Spec A_1^\bullet \underset{\x}{\times} \Spec A_2^\bullet.\]
    Every such pullback is indeed affine because all of our stacks $\x$ are assumed to be separated. As $\F$ is a functor we obtain morphisms of cosimplicial spectra
    \[\F(\Spec A_1^\bullet) \to \F(\Spec B^\bullet) \gets \F(\Spec A_2^\bullet)\]
    respecting any $\E_n$-structure on $\F$.
    The desired comparison morphisms of filtered spectra come from applying the d\'{e}calage functor to the above:
    \[\DSS_1(\x,\F) = \Dec(\F(\Spec A_1^\bullet)) \to \Dec(\F(\Spec B^\bullet)) \gets \Dec(\F(\Spec A_2^\bullet)) = \DSS_2(\x,\F).\]
    To see the left morphism is an equivalence, it suffices to prove it is an isomorphism on bigraded homotopy groups of the associated graded; see \cite[Pr.II.1.9]{alicethesis}, using that these filtrations are complete by \cref{lem:decalage_is_complete}.
    Recall that the bigraded homotopy groups of the associated graded a filtered spectrum are precisely the entries of the first page of its associated spectral sequence.
    By \cref{e2pageofdssidentification}, the homotopy groups of the associated graded of both the source and target of $\Dec(\F(\Spec A_1^\bullet)) \to \Dec(\F(\Spec B^\bullet))$ are naturally isomorphic to the sheaf cohomology of~$\F$.
    The same argument applies to the map $\Dec(\F(\Spec A_2^\bullet)) \to \Dec(\F(\Spec B^\bullet))$.
\end{proof}


\subsection{Descent spectral sequences in synthetic spectra}
\label{ssec:even_periodic_refinement}

In this section, we define the synthetic sheaf $\O^\syn$ advertised in the introduction.
We show that it encodes descent spectral sequences, proving \cref{maingeneral}.
We begin with a discussion of the spectral stacks that our construction applies to: the \emph{even-periodic refinements}.

First, recall that associated to a formal group law $F$ over a ring $R$, the \emph{associated formal group} $\widehat{\mathbf{G}}_F$ represents the functor sending an $R$-algebra to its set of nilpotent elements.
This set becomes a group under the multiplication given by $x+_F y = F(x,y)$.
Using this notation, one can write down the modern interpretation of the Landweber exact functor theorem, as can be found in \cite[Th.0.0.1]{ec2name}, for example.
We let $\M_\FG$ denote the moduli stack of formal groups and let $\omega$ denote the line bundle defined on $\M_\FG$ that sends a formal group $p\colon G\to X$ over a scheme $X$ to the line bundle $p_\ast \Omega^1_{G/X}$.

\begin{theorem}[Landweber exact functor theorem]
Let $R$ be a classical commutative ring with a formal group law $F$; this determines a map of stacks $\Spec R\to \M_\FG$.
If this map is flat, then there exists an even-periodic homotopy commutative ring spectrum $A$ equipped with an isomorphism of rings $\pi_0 A\cong R$ and an isomorphism of formal groups over $R$
\[
    \Spf A^0(\CP^\infty)\cong \widehat{\mathbf{G}}_F.
\]
Moreover, $\pi_{2k}A$ can be naturally identified with the line bundle $\omega^{\otimes k}_R$ associated to ${\widehat{\mathbf{G}}_F}$ over $\Spec R$.
\end{theorem}

The sheaves of $\E_\infty$-rings we study in this article locally look like they came from this theorem.

\begin{mydef}\label{evenperiodicrefinenementdefinition}
    Let $\x$ be a classical Deligne--Mumford stack\footnote{Some of what we do here could be generalised to a larger class of stacks or more general sheaves. We restrict ourselves to Deligne--Mumford stacks here to more easily reference \cite{sagname}. Generalisations are definitely possible.} equipped with a flat map of stacks $f\colon \x\to \M_\FG$.
    Let $\O^\top_\x$ be a sheaf of $\E_\infty$-rings on the small \'{e}tale site of $\x$.
    We say that $(\x,\O^\top_\x)$ is an \emph{even-periodic refinement} (of $\x$ and $f$) if $\O^\top_\x$ is \emph{locally Landweber exact}, meaning that for each \'{e}tale map of stacks $\Spec A\to \x$, the $\E_\infty$-ring $\O^\top_\x(\Spec A)$ represents the Landweber-exact even-periodic cohomology theory associated to the composite
    \[\begin{tikzcd}
        \Spec A\ar[r] & \x \ar[r,"f"] & \M_\FG.
    \end{tikzcd}\]
\end{mydef}

There are many examples of even-periodic refinements.
In this article, we will focus on two main applications.

\begin{example}
    Let $E$ be an even-periodic complex-oriented $\E_\infty$-ring, with an action of a finite group $G$.
    Then the quotient $\Spec (E)/ G$ represents an even-periodic refinement of $\Spec (\pi_0 E)/ G$.
Common examples include the Morava E-theories $\uE_n$ and the action by various finite subgroups of the Morava stabiliser group.
An integral example of this form is $\Spec (\KU)/C_2$ from the introduction.    
\end{example}

\begin{example}
Let $\overline{\M}_\Ell$ be the compactification of the moduli stack of elliptic curves.
Then the Goerss--Hopkins--Miller sheaf $\O^\top$ is an even-periodic refinement of $\overline{\M}_\Ell$ by design; see \cite{bourbakigoerss} and \cite[Th.B]{globaltate}.
More general examples of this flavour include anything coming out of \emph{Lurie's theorem}, see \cite[\textsection5]{luriestheorem}, such as the spectra of topological automorphic forms of Behrens--Lawson \cite{taf}.   
\end{example}

To lift the structure sheaf of an even periodic refinement $\O^\top_\x$ to synthetic spectra, we first show that $\nu$ preserves the sheaf condition on affines.

\begin{prop}\label{nuotopissheaf}
    Let $(\x,\O^\top_\x)$ be an even-periodic refinement.
    Then the composite
    \[
    \begin{tikzcd}
        (\Aff_{/\x}^{\et})^\op \ar[r,"\O^\top_\x"] & \Sp \ar[r,"\nu"] & \Syn_\MU
    \end{tikzcd}
    \]
    is an \'{e}tale sheaf.
\end{prop}

\begin{proof}
    Let $\Spec B^\bullet \to \Spec A$ be a \v{C}ech nerve of an \'{e}tale cover $\Spec B\to \Spec A$. We know $\O^\top_\x$ is an \'{e}tale sheaf, so the natural map of $\E_\infty$-rings
    \[\O^\top_\x(A)\xrightarrow{\cong} \lim \left(   \O^\top_\x(B^0)\Rightarrow \O^\top_\x(B^1)\Rrightarrow \cdots   \right)\]
    is an equivalence, where $\O^\top_\x(\Spec A)=\O^\top_\x(A)$ and likewise for $\O^\top_\x(B^\bullet)$. To see that the natural map of synthetic spectra
    \begin{equation}\label{osynisasheafmap}
        \nu\O^\top_\x(A)\to \lim \left(   \nu\O^\top_\x(B^0)\Rightarrow \nu\O^\top_\x(B^1)\Rrightarrow \cdots   \right)
    \end{equation}
    is an equivalence, we will calculate the bigraded homotopy groups of both sides from the cellularity of $\Syn_\MU$; see \cite[Th.6.2]{syntheticspectra}.
    Since $\O^\top_\x$ on affine schemes is complex orientable, it follows from \cref{prop:nu_of_E_module} that $\pi_{\ast,\ast}\nu\O^\top_\x(A)\cong \pi_\ast \O^\top_\x(A)[\tau]$, and similarly for $\nu\O^\top_\x(B^\bullet)$.
    To calculate the homotopy groups of $\lim \nu \O^\top_\x(B^\bullet)$, we will use a BKSS internal to synthetic spectra. The $\uE_1$-page of this spectral sequence is given by $\pi_\ast \O^\top_\x(B^\bullet)[\tau]$, which is a resolution for the cohomology of $\pi_\ast\O^\top_\x(A)[\tau]$ as $\Z[\tau]$ is flat over $\Z$.
    The $\uE_2$-page of this BKSS is then given by the cohomology of the quasi-coherent sheaf $\pi_\ast\O^\top_\x[\tau]$; this is a quasi-coherent sheaf by \cite[Pr.2.2.6.1]{sagname} or by inspection. However, as $\Spec A$ is affine, the higher cohomology groups of the quasi-coherent sheaf $\pi_\ast \O^\top_\x$ vanish and we see that the edge map
    \[\pi_\ast \O^\top_\x(A)[\tau] \xrightarrow{\cong} \uH^\ast(\Spec A,\, \pi_\ast \O^\top_\x[\tau])\]
    induced by \eqref{osynisasheafmap}, is an isomorphism.
\end{proof}

We can extend this sheaf on affines to a sheaf on all \'etale covers of $\x$; by the above, its value on affine schemes does not change.

\begin{mydef}\label{definitionofsyntheticlift}
    Let $(\x,\O^\top_\x)$ be an even-periodic refinement. We write $\O^\syn_\x$ for the right Kan extension of the sheaf $\opnu \O^\top_\x$ of \cref{nuotopissheaf} to the small \'{e}tale site of $\x$.
    As $\O^\top_\x$ is a sheaf of $\E_\infty$-rings, the sheaf $\O^\syn_\x$ is naturally a sheaf of synthetic $\E_\infty$-rings.
\end{mydef}

Our main object of study is the global sections of $\O^\syn_\x$, denoted by either $\O^\syn_\x(\x)$ or $\Ga(\x,\O^\syn_\x)$. First, we collect some basic properties of this sheaf of synthetic spectra.

\begin{prop}\label{sheafismulocalandtaucomplete}
    The sheaf $\O^\syn_\x$ of \cref{nuotopissheaf} takes values in $\nu \MU$-local {\upshape(}a.k.a.\ hypercomplete{\upshape)} and $\tau$-complete $\MU$-synthetic spectra.
\end{prop}

\begin{proof}
    The classes of $\nu\MU$-local and $\tau$-complete synthetic spectra are both closed under limits, so it suffices to check the statement on affine schemes.
    For every \'etale map $\Spec A\to \x$, the synthetic spectrum $\O^\syn_\x(\Spec A)$ is the synthetic analogue of a homotopy $\MU$-module.
    It is therefore a homotopy $\nu\MU$-module, and in particular $\nu\MU$-local.
    It is $\tau$-complete by \cite[Pr.A.13]{burkhahnseng}.
\end{proof}

Being $\tau$-complete, these synthetic spectra are to some extent fully captured by their signature.
We can compute their signature using the results in \cref{taubsssection}; this gives us \cref{maingeneral}.

\begin{cor}\label{signatureofsyntheticglobalsections}
    Let $(\x,\O^\top_\x)$ be an even-periodic refinement.
    Then there is a equivalence of $\E_\infty$-algebras in filtered spectra
    \[\opsigma \O^\syn_\x(\x) \cong \DSS(\x,\O^\top_\x).\]
    In particular, the signature spectral sequence of $\O^\syn_\x(\x)$ is naturally equivalent to the DSS for $\O^\top_\x(\x)$.
\end{cor}

A more general version of this theorem for quasi-coherent sheaves over an even-periodic refinement appears in \cite{sventhesis}.

\begin{proof}
    This is a special case of \cref{bksssignature}, using \cref{nuotopissheaf}.
\end{proof}

The following is a computationally useful corollary of the above fact. It states that the signature of $\O^\syn_\x(\x)$ always satisfies a \emph{checkerboard} pattern on its $\uE_2$-page.

\begin{cor}
    For an even-periodic refinement $(\x,\O^\top_\x)$, the bigraded homotopy groups of $\O^\syn_\x(\X)$ are \emph{even}, in the sense that their bigraded homotopy groups vanish when the sum of the two bidegrees is odd.
\end{cor}

In \cite{sventhesis}, this is related to even $\MU$-synthetic spectra as in \cite[Df.5.8]{syntheticspectra}.

\begin{proof}
    By \Cref{signatureofsyntheticglobalsections}, the mod $\tau$ homotopy groups of $\O^\syn_\x(\x)$ give the $\uE_2$-page of the DSS for $\O^\top_\x(\x)$. As $\O^\top_\x$ is an even-periodic refinement, this spectral sequence is concentrated in such degrees.
\end{proof}

Furthermore, \cref{prop:limits_are_synthetic_lifts} tells us that $\O^\syn_\x$ is levelwise a synthetic lift of~$\O^\top_\x$.
We can be more precise about their relationship, as follows.
In the following proposition, we use the notation $\nu\circ \O^\top_\x$ to denote the composition of functors; this is not to be confused with the sheafification thereof, which we denote by $\O^\syn_\x$.

\begin{prop}\label{invert_tau_Osyn}
    Let $(\x,\O^\top_\x)$ be an even-periodic refinement.
    \begin{numberenum}
        \item The natural map $\nu\circ \O^\top_\x \to \O^\syn_\x$ factors through an equivalence $\nu\circ \O^\top_\x \cong \tau_{\geq 0} \circ \O^\syn_\x$.
        \item The natural map $\nu\circ \O^\top_\x \to \O^\syn_\x$ becomes an equivalence after applying the functor $\tau^{-1}\colon \Syn_\MU \to \Sp$.
    \end{numberenum}
\end{prop}

\begin{proof}
    This is a special case of \cref{prop:limits_are_synthetic_lifts}, using \cref{nuotopissheaf}.
\end{proof}

In particular, the functor $\nu$ preserves the limit defining $\O^\top_\x(\x)$ if and only if $\O^\syn_\x(\x)$ is connective. 
While this property of being connective may first appear to be a purely abstract condition on a synthetic spectrum, the preceding discussion shows that it has concrete, non-synthetic implications.
Indeed, combining \cref{signatureofgammanu,invert_tau_Osyn} with \cref{signatureofsyntheticglobalsections}, we learn that if $\O^\syn_\x(\x)$ is connective, then there is a natural monoidal equivalence of filtered spectra
\[
    \Dec(\ASS_\MU(\O^\top_\x(\x))) \cong \DSS(\x, \O^\top_\x).
\]
In \cref{sec:comparison_DSS_ANSS}, we determine when $\O^\syn_\x(\x)$ is connective; see \cref{theorem:when_Osyn_connective} and \cref{affine_implies_connective}.


\section{The homology of synthetic global sections}\label{sec:computing_connectivity}

Recall from \cref{invert_tau_Osyn} that the natural map $\nu\O^\top_\x(\x) \to \O^\syn_\x(\x)$ is a connective cover inside $\Syn_{\MU}$.
To check when this map is an equivalence, we need to study the $\nu\MU$-homology of $\O^\syn_\x(\x)$.
In this section, we set up a spectral sequence to compute this homology.
It is an example of a synthetic version of the descent spectral sequence, which we introduce in \cref{sec:sdss}.
Most of the work is in identifying the abutment of this spectral sequence with $\nu\MU\otimes\O^\syn_\x(\x)$; this is the subject of \cref{ssec:identifying_global_sections}.
Once this is set up, we use this spectral sequence to prove \cref{nuandgacommutesometimes} in \cref{sec:comparison_DSS_ANSS}.


\subsection{The synthetic descent spectral sequence}\label{sec:sdss}

The synthetic descent spectral sequence is set up in the exact same way as in \cref{sec:DSS}, but using filtered synthetic spectra and sheafified bigraded homotopy groups instead.

\begin{construction}\label{constr:SDSS}
    Fix a homotopy commutative ring spectrum $E$ of Adams type. Let $\x$ be a Deligne--Mumford stack with $\FF$ be a sheaf of $E$-synthetic spectra on the small \'{e}tale site of~$\x$. For a fixed $E$-synthetic spectrum $S$, let $S \otimes \FF$ denote the sheafification of the presheaf
    \[
        (\DM^\et_{/\x})^\op \to \Syn_{E}, \qquad (\y\to \x) \mapsto S \otimes \FF(\y).
    \]
    The \emph{$E$-synthetic descent spectral sequence} ($E$-SDSS, or SDSS when the role of $E$ is clear) for $S \otimes \FF$ is the DSS for the sheaf $S\otimes \FF$ and takes the form
    \[
        \mathrm{E}_2^{k,w,s} = \mathrm{H}^s(\x,\, \pi_{k+s,\, w-s}(S \otimes \FF)) \implies \oppi_{k,w} \Ga(\x,\, S \otimes \FF),
    \]
    with differentials $d_r^\syn$ of tridegree $(-1,1,r)$.
    Note that multiplication by $\tau$ is a map of tridegree $(0,-1,0)$.
\end{construction}

In the case we care about, the $\uE_2$-page can be made more explicit.

\begin{lemma}
    \label{lem:homology_sheaves_Osyn}
    Let $(\x,\O^\top_\x)$ be an even-periodic refinement and $T$ a spectrum. Then we have an isomorphism of bigraded sheaves of abelian groups on $\x$
    \[
        \pi_{*,*}(\nu T \otimes \O^\syn_\x) \cong \pi_*(T\otimes \O^\top_\x)[\tau],
    \]
    with $\pi_k(T\otimes \O^\top_\x)$ placed in bidegree $(k,0)$.
\end{lemma}

\begin{proof}
    By \cref{prop:nu_of_E_module} and \cref{comeansmumodule}, for every \'etale map $\Spec A \to \x$, we have a natural isomorphism
    \[
        \pi_{*,*}(\nu T\otimes \O^\syn_\x(\Spec A)) \cong \pi_*(T\otimes \O^\top_\x(\Spec A))[\tau].
    \]
    Sheafifying this yields the desired isomorphism.
\end{proof}


\subsection{Identifying the abutment}
\label{ssec:identifying_global_sections}

The goal of this section is to compute the abutment of the SDSS for $\nu\MU\otimes\O^\syn_\x$.
To state the main result, we need some terminology.

\begin{mydef}\label{tamedefinition}
    An even-periodic refinement $(\x,\O_\x)$ is \emph{tame} if the underlying map $\x\to \M_\FG$ is tame; see \cite[Df.2.28]{akhilandlennart} for this relative definition and \cite{tamestacks} for further discussion.
\end{mydef}

\begin{prop}\label{switch_map_is_localization}
    Let $(\x,\O^\top_\x)$ be a tame even-periodic refinement.
    Then the natural limit-comparison map
    \begin{equation}\label{tensorwithmuequation} \nu \MU \otimes  \Ga(\x,\,\O^\syn_\x) \to  \Gamma(\x,\, \nu \MU \otimes \O^\syn_\x)\end{equation}
    is an equivalence.
\end{prop}

\begin{remark}
\label{rmk:no_synthetic_affineness}
The previous proposition describes an instance when the synthetic global sections functor preserves a colimit. 
Mathew--Meier study the spectral global sections functor
\[
    \Gamma^\top \colon \QCoh(\x,\O_\x^\top)\to\Mod_{\Gamma(\x,\O_\x^\top)}(\Sp).
\]
One of their main results  \cite[Th.4.14]{akhilandlennart} states that if $\x$ is tame, then the functor $\Gamma^\top$ preserves \textit{all} colimits.
They moreover give a sufficient condition for $\Gamma^\top$ to be an equivalence; see \cite[Th.4.1]{akhilandlennart}.

These affineness results fail synthetically, even in the simplest examples.
Consider, as in the introduction, the even-periodic refinement $\Spec(\KU)/C_2$ of the affine map $BC_2\to\M_\FG$.
The corresponding $\infty$-category of quasicoherent synthetic sheaves may be identified with $\Mod_{\opnu\KU}(\Syn_\MU)^{hC_2}$, as in 
\cite[Pr.2.16]{akhilandlennart}.
We claim that the global sections functor
\[\Gamma\colon \Mod_{\opnu\KU}(\Syn_\MU)^{hC_2}\to \Mod_{\nu\KO}(\Syn_\MU)\]
does not preserve all colimits, as can be seen by the non-nilpotence of $\eta\in\oppi_{1,1}\nu\KO$.
Indeed, the unit $\O^\syn\in \Mod_{\nu\KU}(\Syn_\MU)^{hC_2}$ has a class $\eta\in\oppi_{1,1}\O^\syn$ lifting $\eta\in\oppi_{1,1}\nu\KO$ along the global sections functor, by linearity over the synthetic sphere $\Sph\in\Syn_\MU$.
The forgetful functor
\[\Mod_{\nu\KU}(\Syn_\MU)^{hC_2}\to \Mod_{\nu\KU}(\Syn_\MU)\]
is conservative, preserves colimits, and sends $\eta$ to $0$, hence
\[\O^\syn[\eta^{-1}] = 0 \quad \text{in } \Mod_{\nu\KU}(\Syn_\MU)^{hC_2}.\]
The colimit-comparison map
\[\Gamma(\O^\syn)[\eta^{-1}]=(\nu\KO)[\eta^{-1}]\to \Gamma(\O^\syn[\eta^{-1}])\]
therefore cannot be an equivalence, since $\eta\in\oppi_{1,1}\nu\KO$ is not nilpotent.
\end{remark}

\begin{remark}\label{tmfcounterexample}
A similar analysis shows that the synthetic global sections functor
\[\Ga\colon \QCoh(\M_\Ell, \O^\syn) \to \Mod_{\O^\syn(\M_\Ell)}(\Syn_{\MU})\]
is \textbf{not} an equivalence.
Again, the reason is that $\oppi_{1,1}\O^\syn(\M_\Ell)$ contains the nonnilpotent element $\eta$.
Indeed, in \cref{sec:SMF}, it is shown that $\O^\syn(\M_\Ell)\cong \opnu\TMF$, from which it follows that this synthetic spectrum contains the nonnilpotent element $\eta$.
\end{remark}

The key problem occurring in the above remarks is one of \emph{hypercompleteness}: if we consider the synthetic global sections functor taking values in the $\infty$-category of modules over $\O^\syn_\X(\X)$ in hypercomplete $\MU$-synthetic spectra, then one does obtain various cocontinuity statements; see \cite{sventhesis} for more.

To prove \cref{switch_map_is_localization}, we will check that it is an equivalence after inverting (\cref{tensorwithmugenericfibre}) and killing (\cref{cor:equiv_after_modtau}) the distinguished element $\tau$.
The first case reduces to the spectral global sections considered by Mathew--Meier.



\begin{lemma}\label{tensorwithmugenericfibre}
    The map~\eqref{tensorwithmuequation} becomes, after $\tau$-inversion, the limit-comparison map
    \[\MU \otimes \Gamma^\top(\x,\O^\top_\x) \to \Gamma^\top(\x,\, \MU \otimes \O^\top_\x).\]
    Moreover, if $(\x,\O^\top_\x)$ is tame, then the above map is an equivalence.
\end{lemma}

\begin{proof}
    The first claim follows from \cref{prop:limits_are_synthetic_lifts}, using that the functor $\tau^{-1} \colon \Syn \to \Sp$ is symmetric monoidal.
    The second claim follows from \cite[Th.4.14]{akhilandlennart}.
\end{proof}

The case of killing $\tau$ is more delicate, as the synthetic global sections do not preserve all colimits by \cref{rmk:no_synthetic_affineness,tmfcounterexample}.
Nevertheless, our proof for the mod $\tau$ version of \cref{switch_map_is_localization} is inspired by the proofs of \cite[Prs.3.9, 4.11]{akhilandlennart}.\\


We begin by identifying the map obtained by tensoring \eqref{tensorwithmuequation} with $C\tau$.
We will consider tensoring with $C\tau$ as a functor
\[-\otimes C\tau\colon \Syn\to \Mod_{C\tau}(\Syn) \cong \Stable_{\MU_\ast\MU},\]
where the equivalence comes from \cite[Pr.4.53]{syntheticspectra}. 
Let us introduce some notation.

\begin{mydef}
Write $\calG$ for the functor
\[
    \calG \colon \Sp \to \Comod_{\MU_*\MU}, \quad X \mapsto \MU_*(X).
\]
We can consider $\calG(X)$ as an object of $\Stable_{\MU_*\MU}$, and we will do so without altering the notation.
If $\O$ is a sheaf of spectra on a stack $\x$, then we denote by $\calG(\O)$ the sheafification in $\Stable_{\MU_*\MU}$ of the composite
\[
    \begin{tikzcd}
        (\DM_{/\x}^\et)^\op \ar[r,"\O"] & \Sp \ar[r,"\calG"] & \Stable_{\MU_*\MU}.
    \end{tikzcd}
\]
\end{mydef}

\begin{lemma}
    \label{lem:identifcation_comp_map_after_Ctau}
    Let $(\x,\O^\top_\x)$ be an even-periodic refinement.
    Then the functors $\calG(\O^\top_\x)$ and $\O^\syn_\x \otimes C\tau$ are naturally equivalent.
    Moreover, the map~\eqref{tensorwithmuequation} after tensoring with $C\tau$ can be identified {\upshape(}as a map in $\Stable_{\MU_*\MU}${\upshape)} with the limit-comparison map
    \[
        \calG(\MU) \otimes \Gamma(\x,\, \calG(\O^\top_\x)) \to \Gamma(\x,\, \calG(\MU) \otimes \calG(\O^\top_\x)).
    \]
\end{lemma}

\begin{proof}
    Note that the functor
    \[-\otimes C\tau:\Syn\to\Mod_{C\tau}(\Syn) \cong \mathrm{Stable}_{\MU_*\MU}\]
    is symmetric monoidal and preserves limits (as $C\tau$ is dualizable), and that by \cite[Rmk.4.55]{syntheticspectra} we have a natural equivalence $\nu({-})\otimes C\tau \cong \calG$.
    Hence $\O^\syn_\x \otimes C\tau$ is a sheaf, and we can identify it with $\calG(\O^\top_\x)$ because the two sheaves agree on affines.
    This then implies the second claim.
\end{proof}

To prove that \eqref{tensorwithmuequation} induces an equivalence mod $\tau$, we will actually work with the subcategory $\D(\Comod_{\MU_\ast\MU})$ of $\Stable_{\MU_\ast \MU}$, the full subcategory of $\nu\MU$-local objects; see \cite[Th.4.54]{syntheticspectra}.
This is sufficient as both sides of \eqref{tensorwithmuequation} are $\nu\MU$-local (the first as it is a $\nu\MU$-module, the second by \cref{sheafismulocalandtaucomplete}).
To avoid confusion, let us introduce some notation.

\begin{notation}
    For $M,N \in \D(\Comod_{\MU_*\MU})$, write
    \[
        M\hotimes N \defeq L_{\nu\MU}(M\otimes N)
    \]
    for the $\nu\MU$-local tensor product, formed in $\Stable_{\MU_*\MU}$.
    This is again an object of the subcategory $\D(\Comod_{\MU_*\MU})$; in fact, this results in a symmetric monoidal structure on the derived  $\infty$-category.
\end{notation}

Let $f \colon \x \to \M_\FG$ be the map underlying our even-periodic refinement, let $L$ denote the Lazard ring, and let $q \colon \Spec L \to \M_\FG$ denote the cover.
We form the pullback
\begin{equation}\label{pullbackwithspecl}
    \begin{tikzcd}
        \y \ar[r,"\bar{f}"] \ar[d,"\bar{q}"'] \pullback & \Spec L \ar[d,"q"] \\
        \x \ar[r,"f"] & \M_\FG.
    \end{tikzcd}
\end{equation}
We remind the reader that $q$ is a faithfully flat and affine map of stacks.
The map $\bar{q}$ is therefore faithfully flat and affine as well. In the following statement, the pullback and pushforward functors denote their derived variants, landing in the derived $\infty$-category of comodules.\\

For the following lemma, let $\omega$ be the line bundle on $\M_\FG$, defined to send some $\Spec R\to \M_\FG$ defining a (1-dimensional) formal group $\mathbf{G}$ over $R$ to the dual of $\ker(\mathbf{G}(R[t]/t^2)\to \mathbf{G}(R))$, i.e., the dual to the Lie algebra of $\mathbf{G}$. We then write $\O_{\Spec L}$ for the sheaf of graded abelian groups given by the pullback $q^\ast\omega^{\otimes\ast}$.

\begin{lemma}\label{lemma:beckchevalley}
    With $\x$ as above, there are equivalences of quasi-coherent sheaves on the classical Deligne--Mumford stack $\x$
    \[
        \calG(\O^\top_\x) \cong f^\ast q_\ast \O_{\Spec L} \cong \bar{q}_* \bar{f}^*\O_{\Spec L} \cong \bar{q}_\ast \bar{q}^\ast \omega_\x^{\otimes \ast/2}.
    \]
    In particular, for every $s\geq 0$ and every $t$, we have isomorphisms
    \[
        \uH^s(\x,\, \pi_t(\MU\otimes\O^\top_\x)) \cong \uH^s(\y,\, \bar{q}^\ast \omega^{\otimes t/2}_\x).
    \]
\end{lemma}

\begin{proof}
    For the first equivalence, both sides are sheaves, so it suffices to show they agree on affines. As $f$ is flat and $q$ is flat and affine, $f^*$ and $q_*$ are exact functors, so the derived sheaf $f^*q_*\O_{\Spec L}$ sends an affine $\Spec A_0\to \x$ to the ($2$-periodic) $\MU_*\MU$-comodule algebra appearing as the pullback
    \[
    \Spec A_0\times_{\M_\FG}\Spec L 
    \]
    with its comodule structure from the pullback descent datum for the cover $q$, regarded as a cochain complex concentrated in cohomological degree zero. This pullback is computed as $\MU_*A$ by \cite[Pr.2.4]{akhilandlennart} where $A$ is the Landweber theory associated to the flat map $\Spec A_0\to\M_\FG$.\\

    For the second equivalence, we refer to \cite[Cor.2.5.4.6]{sagname}, which says that the canonical Beck--Chevalley transformation $f^\ast q_\ast\to \bar{q}_\ast \bar{f}^\ast$ is an isomorphism. The third equivalence is simply the observations that the sheaves $\pi_t \O^\top_\x \cong \omega_\x^{\otimes t/2}$ are isomorphic by design, and the consequence that $\pi_\ast\bar{q}^\ast \O^\top_\x\cong \pi_\ast\O^\top_\y \cong \bar{f}^\ast \O_{\Spec L}$, as all of these sheaves pull back from $\omega^{\otimes \ast}$ over $\M_\FG$. The isomorphism on cohomology follows as $\bar{q}$ is affine and flat; see, e.g., \cite[Lm.6.2.10]{carrickthesis}.
\end{proof}

\begin{mydef}\label{qcohingendef}
    Let $\C$ be a presentably symmetric monoidal $\infty$-category, $\x$ a classical Deligne--Mumford stack, and $\O$ a sheaf on the small \'{e}tale site of $\x$ with values in $\CAlg(\C)$.
    We say that a module $M$ over $\O$ is \emph{quasi-coherent} if for every map of affine schemes $U\to V$ in the small \'{e}tale site of $\x$, the natural map of $\O(U)$-modules
    \[M(V)\underset{\O(V)}{\otimes} \O(U)\to M(U)\]
    is an equivalence. Let us write $\QCoh(\x,\O)$ for the full $\infty$-subcategory of $\Mod_\O(\Shv_\C(\x))$ spanned by the quasi-coherent sheaves.
\end{mydef}

In the next lemma, we consider $\calG(\O^\top_\x)$ as a sheaf valued in $\D(\Comod_{\MU_*\MU})$.

\begin{lemma}
    \label{lem:glob_sect_preserve_colims}
    Let $(\x,\O^\top_\x)$ be as above.
    If $(\x,\O^\top_\x)$ is tame, then the functor
    \[
        \QCoh(\x,\calG(\O^\top_\x)) \to \D(\Comod_{\MU_*\MU}), \qquad F \mapsto \Gamma(\x,\, F \hotimes \calG(\O^\top_\x))
    \]
    preserves small colimits.
\end{lemma}


\begin{proof}
    The functor is an exact functor between stable $\infty$-categories, so it preserves finite colimits.
    It suffices therefore to show that it preserves small coproducts. We follow closely the argument of \cite[Pr.3.9]{akhilandlennart}, and we make use of the SDSS
    \begin{equation}\label{equation:descentsscoproduct}
    \uE_2^{k,w,s} =\uH^s(\x,\, \pi_{k+s,\, w-s}(F\hotimes \calG(\O^\top_\x)))\implies \oppi_{k,w}\Gamma(\x,\, F\hotimes\calG(\O^\top_\x)).
    \end{equation}
The $\uE_2$-page of this SS commutes with coproducts in $F$ since $\x$ is quasi-compact and separated (see \cite[Lm.3.10]{akhilandlennart}), hence so does the $\uE_\infty$-page. If the $\uE_\infty$-page is concentrated in finitely many rows, then the filtrations are finite, and it follows that the natural map
\[\bigoplus\limits_i \Gamma(\x,\, F_i\hotimes\calG(\O^\top_\x))\to\Gamma(\x,\, \bigoplus\limits_i F_i\hotimes\calG(\O^\top_\x))\]
induces an isomorphism on bigraded homotopy groups and is therefore an equivalence.\\

It suffices to show that the $\uE_2$-page of \eqref{equation:descentsscoproduct} is concentrated in finitely many rows $s<N$ for fixed $N$ independent of $F$. This follows from the chain of isomorphisms
\begin{align*}
    \uH^s(\x,\, \pi_{*,*}(F\hotimes\calG(\O^\top_\x)))&\cong \uH^s(\x,\, \pi_{*,*}(F\hotimes \bar{q}_* \bar{f}^*\O_{\Spec L}))\\
    &\cong \uH^s(\x,\, \pi_{*,*}(\bar{q}_*(\bar{q}^*F\hotimes \bar{f}^*\O_{\Spec L})))\\
    &\cong \uH^s(\x,\, \bar{q}_*\pi_{*,*}(\bar{q}^*F\hotimes \bar{f}^*\O_{\Spec L}))\\
    &\cong \uH^s(\y,\, \pi_{*,*}(\bar{q}^*F\hotimes \bar{f}^*\O_{\Spec L})).
\end{align*}
The first isomorphism is \cref{lemma:beckchevalley}, the second is the projection formula, and the third and fourth follow since $\bar{q}$ is affine and flat (see \cite[Lm.6.2.10]{carrickthesis} for example). The claim now follows from the fact that $\y$ has finite cohomological dimension.
\end{proof}

\begin{cor}
    \label{cor:equiv_after_modtau}
    Let $(\x,\O^\top_\x)$ be a tame even-periodic periodic refinement. Then the map \eqref{tensorwithmuequation} becomes an equivalence after tensoring with $C\tau$.
\end{cor}
\begin{proof}
    Let $\C$ denote the full subcategory of $\D(\Comod_{\MU_*\MU})$ on those $M$ such that the limit-comparison map
    \begin{equation}
        M \hotimes \Gamma(\x,\, \calG(\O^\top_\x)) \to \Gamma(\x,\, M \hotimes\calG(\O^\top_\x))\label{eq:swapping_Gamma_with_tensor}
    \end{equation}
    is an equivalence.
    Clearly $\C$ contains all dualisable objects; in particular, it contains the bigraded spheres $\Sph^{k,s}\otimes C\tau$.
    Moreover, we claim that $\C$ is closed under small ($\nu\MU$-local) colimits.
    Indeed, suppose that we are given a diagram $M_\alpha$ in $\C$. Then for $M= \colim_\al M_\al$, the map~\eqref{eq:swapping_Gamma_with_tensor} may be factored as
    \[
        \colim_\alpha M_\alpha \hotimes \Gamma(\x, \, \calG(\O^\top_\x)) \to \colim_\alpha \Gamma(\x, \, M_\alpha \hotimes \calG(\O^\top_\x)) \to \Gamma(\x, \,\colim_\alpha M_\alpha \hotimes \calG(\O^\top_\x)).
    \]
    The first of these maps is an equivalence, being a colimit of equivalences, and the second is an equivalence as a consequence of \cref{lem:glob_sect_preserve_colims}.
    
    The $\infty$-category $\D(\Comod_{\MU_*\MU})$ is cellular in the sense that it is generated under ($\nu\MU$-local) colimits by the bigraded spheres.
    Hence $\C$ is equal to all of $\D(\Comod_{\MU_*\MU})$.
    In particular, the map \eqref{eq:swapping_Gamma_with_tensor} for $M = \calG(\MU)$ is an equivalence.
    For this choice of $M$, the $\nu\MU$-local tensor products coincide with the ordinary ones formed in $\Stable_{\MU_*\MU}$.
    Therefore, by \cref{lem:identifcation_comp_map_after_Ctau}, the map \eqref{eq:swapping_Gamma_with_tensor} for $M = \calG(\MU)$ coincides with the map \eqref{tensorwithmuequation} after tensoring with $C\tau$, showing that the latter is an equivalence.
\end{proof}

\begin{proof}[Proof of \cref{switch_map_is_localization}]
    The desired map \eqref{tensorwithmuequation} is an equivalence after inverting $\tau$ and tensoring with $C\tau$ by \cref{tensorwithmugenericfibre} and \cref{cor:equiv_after_modtau}, respectively, which implies that it is an equivalence.
\end{proof}

\begin{remark}
    The same proof of \cref{switch_map_is_localization} would work with $\nu\MU$ replaced by a homotopy $\nu\MU$-module $M$ in synthetic spectra.
    Indeed, this condition ensures that $M/\tau \hotimes \calG(\O^\top_\x) \cong M/\tau \otimes \calG(\O^\top_\x)$ as objects of $\Stable_{\MU_*\MU}$.
    For example, this applies when $M$ is the synthetic analogue of a complex-orientable ring spectrum.
    For a general synthetic spectrum $M$ in the place of $\nu\MU$, our proof shows that the map \eqref{tensorwithmuequation} is $\nu\MU$-localisation.
\end{remark}


\subsection{Comparing the DSS and the ANSS}
\label{sec:comparison_DSS_ANSS}

To use the SDSS for $\nu\MU\otimes \O^\syn_\x$ to prove \cref{nuandgacommutesometimes}, we will need more precise information about this spectral sequence. The following is in the vein of \cite[Lm.15]{lemmafifteen} and \cite[Th.A.8]{burkhahnseng}.

\begin{prop}
    \label{thm:SDSS}
    Let $(\x,\O^\top_\x)$ be a tame even-periodic refinement.
    Write $\{\mathrm{E}_r^{k,w,s},d^\syn_r\}_r$ for the SDSS for $\nu \MU \otimes \O^\syn_\x$, and write $\{\mathrm{E}_r^{k,s},d^\top_r\}_r$ for the DSS for $\MU \otimes \O^\top_\x$.
    \begin{numberenum}
        \item The SDSS for $\nu \MU \otimes \O^\syn_\x$ converges strongly to $\nu \MU_{k,w}(\O^\syn_\x(\x))$.\label{claim:converge}
        \item Using the notation of \eqref{pullbackwithspecl}, we have an isomorphism of trigraded $\Z[\tau]$-modules\footnote{As usual, we set $\omega^{n/2}=0$ for odd $n$.}
        \[
            \mathrm{E}_2^{k,*,s} \cong \mathrm{E}_2^{k,s} \otimes\Z[\tau] \cong \mathrm{H}^s(\y, \omega_\y^{\otimes(k+s)/2}) \otimes \Z[\tau],
        \]
        where $\mathrm{E}_2^{k,s}$ is placed in tridegree $(k,s,s)$, and $\tau$ has tridegree $(0,-1,0)$.\label{claim:E2_SDSS}
        \item The differentials in the SDSS are determined by those in the DSS for $\MU \otimes \O^\top_\x$ in the following way.
        Given a differential $d^\top_r(x) = y$ in the DSS, there is a differential $d^\syn_r(x) = \tau^{r-1}y$ in the SDSS, and all differentials arise in this way.\label{claim:diffs_in_SDSS}
    \end{numberenum}
\end{prop}

\begin{proof}
    Part~\ref{claim:E2_SDSS} is the combination of \cref{lem:homology_sheaves_Osyn} with \cref{lemma:beckchevalley}.
    Part~\ref{claim:diffs_in_SDSS} is proved in the exact same way as \cite[Th.A.8]{burkhahnseng}.
    This leaves us with the claim about convergence.
    The spectral sequence converges conditionally to the global sections of $\nu \MU \otimes \O^\syn_\x$, which by \cref{switch_map_is_localization} is equivalent to $\nu \MU \otimes \Gamma(\x,\O^\syn_\x)$. This gives us conditional convergence to the desired abutment. Moreover, because $\x \to \M_\FG$ is tame, the stack $\y$ has finite cohomological dimension by \cite[Pr.2.29]{akhilandlennart}.
    This means that the SDSS has a horizontal vanishing line (in the $s$ variable) on the $\uE_2$-page, so that the convergence is in fact strong.
    This gives us part~\ref{claim:converge}.
\end{proof}

The connection to the algebraic geometry of $\x$ allows us to answer the question of the connectivity of $\O^\syn_\x(\x)$ in purely algebro-geometric terms.
More specifically, the fact that $\y$ has finite cohomological dimension implies that $\O^\syn_\x(\x)$ is bounded below, and we can in fact determine the precise bound, as follows.

\begin{theorem}
    \label{theorem:when_Osyn_connective}
    Let $(\x,\O^\top_\x)$ be a tame even-periodic refinement, and let $\y$ denote the pullback as in \eqref{pullbackwithspecl}.
    Then for every $n\geq 0$, the following are equivalent:
    \begin{letterenum}
        \item the $\MU$-synthetic spectrum $\O^\syn_\x(\x)$ is $(-n)$-connective;\label{cond:Osyn_connective}
        \item the sheaf cohomology groups $\mathrm{H}^s(\y, \omega_\y^{\otimes t})$ vanish for all $t\in\Z$ and all $s>n$.\label{cond:sheaf_cohom_vanishes}
    \end{letterenum}
    In particular, $\O^\syn_\x(\x)$ is $(-N)$-connective, where $N$ is the cohomological dimension of~$\y$.
\end{theorem}

\begin{proof}
    Recall that an $\MU$-synthetic spectrum $S$ is $(-n)$-connective if and only if $\nu\MU_{\ast,s}(S)$ vanishes for all $s > n$.
    We compute the $\nu\MU$-homology of $\O^\syn_\x(\x)$ using the SDSS from \cref{thm:SDSS}.
    As a result, we learn that $\O^\syn_\x(\x)$ is $(-n)$-connective if and only if the $\mathrm{E}_\infty$-page of the SDSS vanishes in those tridegrees $(k,w,s)$ where $w > n$.
    
    Using part~\ref{claim:E2_SDSS} of \cref{thm:SDSS}, we see that assumption~\ref{cond:sheaf_cohom_vanishes} is equivalent to the analogous vanishing condition for the $\uE_2$-page, so \ref{cond:sheaf_cohom_vanishes} certainly implies \ref{cond:Osyn_connective}.
    To prove the converse, suppose that the group $\mathrm{H}^s(\y,\omega_\y^{\otimes t})$ does not vanish for certain $s,t$, where $s>n$.
    By part~\ref{claim:E2_SDSS} of \cref{thm:SDSS}, this group appears on the $\mathrm{E}_2$-page of the SDSS in tridegree $(2t-s,\, s,\, s)$.
    For degree reasons it cannot be the target of a differential, as differentials only hit elements that are $\tau$-divisible on the $\uE_2$-page, and this tridegree is not $\tau$-divisible.
    If this tridegree contains a permanent cycle, then this gives a nonzero contribution to stem $(2t-s,\, s)$, showing that $\O^\syn_\X(\X)$ is not $(-n)$-connective since $s>n$.
    If, on the other hand, all elements in this tridegree support differentials, then let us pick an element $x$ that supports a $d_r^\syn$-differential for some $r\geq 2$.
    Then $d_r^\syn(x)=\tau^{r-1}\cdot y$ for some class $y$ in tridegree $(2t-s-1,\,s+r,\,s+r)$.
    By $\tau$-linearity of the differential $d_r^\syn$, it follows that $y$ is a permanent cycle and gives a nonzero contribution to stem $(2t-s-1,\, s+r)$.
    Since $s+r>n$, it follows that $\O^\syn_\X(\X)$ is not $(-n)$-connective.
\end{proof}

There are some practical situations when the conditions of \cref{theorem:when_Osyn_connective} are satisfied.

\begin{cor}\label{affine_implies_connective}
    Let $(\x,\O^\top_\x)$ be an even-periodic refinement such that the underlying map $\x\to \M_\FG$ is affine. Then the natural map of synthetic $\E_\infty$-rings
    \[\nu (\O_\x^\top(\x)) \to \O^\syn_\x(\x)\]
    is an equivalence.
    In particular, this induces an equivalence of $\E_\infty$-algebras in filtered spectra
    \[\Dec (\ASS_{\MU}(\Ga(\x,\O^\top_\x))) \cong \DSS(\x,\O^\top_\x).\]
\end{cor}

\begin{proof}
    As the map $\x\to \M_\FG$ is affine, the fibre product $\y =\x \times_{\M_\FG}\Spec L$ is also affine.
    In particular, $\y$ has cohomological dimension zero.
    Affine maps are tame, so \cref{theorem:when_Osyn_connective} applies, showing that the synthetic $\E_\infty$-ring $\O^\syn_\x(\x)$ is connective.
    By \cref{invert_tau_Osyn}, we see that the natural map of synthetic $\E_\infty$-rings is an equivalence.
    Applying $\sigma\colon \CAlg(\Syn_{\MU}) \to \CAlg(\Fil(\Sp))$ to this map, we obtain the middle equivalence of $\E_\infty$-objects in filtered spectra
    \[\Dec (\ASS_{\MU}(\Ga(\x,\O^\top_\x))) \cong \sigma \nu (\O_\x^\top(\x)) \cong \opsigma \O^\syn_\x(\x) \cong \DSS(\x,\O^\top_\x),\]
    the other two equivalences following from \cref{signatureofgammanu} and \cref{signatureofsyntheticglobalsections}, respectively.
\end{proof}

This identification of spectral sequences, and other similar statements, have appeared elsewhere in the literature:

\begin{itemize}
    \item In \cite[Cor.5.3]{tmfhomology}, Mathew shows that the ANSS for $\tmf$ takes the form of a descent spectral sequence (or rather, its $\uE_2$-page is given by the cohomology of a particular Hopf algebroid).
    \item In unpublished works, Meier \cite[Th.4.6]{relativelyfreetmf} and Devalapurkar \cite[Cor.5]{equivwoods} prove the above statement as an isomorphism of spectral sequences, but without addressing the multiplicative structures.
\end{itemize}


\section{Applications}\label{sec:application}

We return now to the examples discussed in the introduction.


\subsection{Homotopy fixed-point spectral sequences}\label{sec:HFPSS}

Let $E$ be an even-periodic Landweber-exact $\E_\infty$-ring with an action by a finite group $G$, that is, $E$ is a functor $BG \to \CAlg(\Sp)$. Then there is a canonical morphism of stacks
\[\Spec(\pi_0E)/G\to\M_\FG\]
that admits a canonical even-periodic refinement given by $\Spec (E)/G$; see \cite[Pr.2.15]{akhilandlennart}.
It is well known that the homotopy fixed-point spectral sequence (HFPSS) for $E$ can be naturally identified with the DSS for the global sections of this even-periodic refinement, which we write as $(\x,\O^\top_\x)$ for a moment.
Indeed, the cobar resolution computing the homotopy fixed-points $E^{hG}$ takes the form of the cosimplicial object
\begin{equation}\label{cobarforhfpss} E^\bullet = \left(E\Rightarrow \prod_G E \Rrightarrow \prod_{G^2} E \cdots\right),\end{equation}
so $E^{hG}$ is given by $\lim E^\bullet$. Notice that we have $E=\O^\top_\x(\Spec \pi_0 E)$ by definition, and that there are identifications
\[
    \O^\top_\x\left(\coprod_{G^n} \Spec \pi_0 E\right) \cong \prod_{G^n} E \quad\text{and}\quad \coprod_{G} \Spec \pi_0 E \cong \Spec \pi_0 E \underset{\Spec \pi_0 E/G}{\times} \Spec \pi_0 E,
\]
the second coming from the fact that $\Spec \pi_0 E\to \Spec (\pi_0 E)/G$ is an \'{e}tale $G$-torsor.
Inductively, we now see that \eqref{cobarforhfpss} is naturally equivalent to the cosimplicial $\E_\infty$-ring $\O^\top_\x(\Spec (\pi_0 E)^\bullet)$ used to define the DSS for the global sections of $(\x,\O^\top_\x)$, where $(\pi_0 E)^\bullet$ comes from the \v{C}ech nerve of above $G$-torsor.\\

In particular, combining this discussion with \cref{comeansmumodule} and \cref{bksssignature} we obtain the following corollary.

\begin{cor}\label{hfpsssignature}
    Let $G$ be a finite group and let $E\in \mathrm{Fun}(BG,\mathrm{CAlg}(\Sp))$ whose underlying ring spectrum is even-periodic and Landweber exact.
    Then, using the notation above, there is a natural equivalence of cosimplicial $\E_\infty$-rings $E^\bullet \cong \O^\top_\x(\Spec (\pi_0 E)^\bullet)$.
    In particular, the HFPSS for $E^{hG}$ and the DSS for $\x$ are isomorphic, and the signature spectral sequence of the synthetic $\E_\infty$-ring $(\nu E)^{hG}$ is the HFPSS for $E$.
\end{cor}

Our framework gives a straightforward criterion that allows one to often identify this DSS (and hence also the HFPSS) with the ANSS of $E^{hG}$.

\begin{theorem}\label{thm:hfpss}
    Let $G$ be a finite group and let $E\in \mathrm{Fun}(BG,\mathrm{CAlg}(\Sp))$ whose underlying ring spectrum is even-periodic and Landweber exact.
    If, for every field-valued point $x$ of $\Spec \pi_0E$, the subgroup of $G$ stabilising $x$ acts faithfully on the formal group associated to $x$, then the canonical map
    \[\nu(E^{hG})\to(\nu E)^{hG}\]
    is an equivalence of synthetic $\E_\infty$-rings. In particular, the signature spectral sequence of $(\nu E)^{hG}$ is the ANSS of~$E^{hG}$.
\end{theorem}

\begin{proof}
    As remarked in \cite[\textsection 6.2]{akhilandlennart}, the morphism 
    \[\Spec(\pi_0E)/G\to\M_\FG\]
    is affine precisely when the conditions of the theorem hold. The claim now follows immediately from \cref{nuandgacommutesometimes}.
\end{proof}

\begin{remark}
    In both \Cref{hfpsssignature} and \Cref{thm:hfpss}, it is not important that $E$ is even-periodic and Landweber exact, only that it is a complex-orientable spectrum.
\end{remark}

This theorem is a highly structured lift of the following useful corollary, which we now see is immediate from \cref{nuandgacommutesometimes}.

\begin{cor}
    Under the conditions of \cref{thm:hfpss}, the ANSS of $E^{hG}$ is isomorphic to the HFPSS of $E$ from the $\uE_2$-page on.
\end{cor}

\begin{example}\label{example:KObaby}
    As the $C_2$-acting on $\KU$ is a subgroup of the automorphism group of $\widehat{\G}_m$, it follows that $G$ always acts faithfully on geometric points of $\Spec \Z$. In particular, it follows that the ANSS for $\KO$ is isomorphic to the $C_2$-homotopy fixed-point spectral sequence of $\KU$, a well-known fact. The equivalence of synthetic $\E_\infty$-rings $\nu \KO \to (\nu\KU)^{h C_2}$ was also shown empirically in \cite[Pr.4.25]{syntheticj}.
\end{example}

\begin{example}\label{example:eon}
    Fix a height $n$ formal group $\Gamma$ over a perfect field $k$ of characteristic $p$, and fix a finite subgroup $G$ of the corresponding Morava stabilizer group $\mathbf{G}\defeq\mathrm{Gal}(k/\mathbf{F}_p)\rtimes\mathrm{Aut}(\Gamma)$. By the Goerss--Hopkins--Miller theorem, see \cite[\textsection5]{ec2name}, the even-periodic Landweber-exact $\E_\infty$-ring $\uE_n$ admits an action of $G$ by $\E_\infty$-ring maps. The conditions of \cref{thm:hfpss} apply, and we get an equivalence of synthetic $\E_\infty$-rings
    \[\nu(\EO_n(G))\cong\nu(\uE_n)^{hG},\]
    where $\EO_n(G)\defeq\uE_n^{hG}$ is a higher real K-theory.
\end{example}

\begin{remark}
    In general, regardless of whether or not $\Spec (\pi_0 E)/G \to \M_\FG$ is affine, the signature spectral sequence of the synthetic spectrum $\nu(E)^{hG}$ is the HFPSS of $E$ by \cref{hfpsssignature}.
    This will not generally coincide with the ANSS of $E^{hG}$, making $\nu(E)^{hG}$ often preferable over $\nu(E^{hG})$ as a synthetic lift of $E^{hG}$.
    The unit of the synthetic $\E_\infty$-ring $\nu(E)^{hG}$ still provides a map of multiplicative spectral sequences
    \[\mathrm{ANSS}(\Sph)\to\mathrm{HFPSS}(E)\]
    suitable for detection arguments, for example. In this sense, the HFPSS for $E$ is a \emph{modified ANSS} \`{a} la \cite[Df.2.6]{syntheticj}. From this point of view, the ANSS of $E^{hG}$ has no advantage over the HFPSS.
\end{remark}

\begin{example}
    Consider $\KU$ with a trivial $C_2$-action. The spectral Deligne--Mumford stack $\Spec(\KU)/C_2$ is an even-periodic refinement of the morphism 
    \[f\colon BC_2\to\M_\FG\]
    that picks out $\widehat{\G}_m$ with a trivial $C_2$-action. The morphism $f$ is not faithful, and hence does not satisfy the conditions of \cref{thm:hfpss}. 
    
    In fact, the corresponding homotopy fixed-point spectral sequence does not coincide with the Adams--Novikov spectral sequence of $\KU^{hC_2}\cong \KU^{B{C_2}_+}$. Indeed, the latter is complex-oriented as a $\KU$-module, and hence its ANSS is concentrated in filtration zero. On the other hand, the cohomology groups $\uH^s(C_2,\pi_*\KU)$ are nontrivial for positive values of $s$, hence the homotopy fixed-point spectral sequence is nontrivial in positive filtrations.

    Since the ANSS of $\KU^{B{C_2}_+}$ collapses on the zero-line, the ANSS provides no computational leverage; computing the $\uE_2$-page already requires knowing $\pi_*(\KU^{B{C_2}_+})$. The HFPSS, however, has an $\uE_2$-page that is algebraically computable using only $\pi_*\KU$ as input. 
\end{example}


\subsection{Synthetic modular forms}
\label{sec:SMF}

Finally, we would like to apply the techniques of this paper to the theory of \emph{topological modular forms}. These objects motivated much of this article, and provide both examples and counterexamples to \cref{nuandgacommutesometimes}.\\

Recall from \cite[\textsection12]{tmfbook} that there exists an even-periodic refinement $\O^\top$ on the \emph{compactification of the moduli stack of elliptic curves} $\overline{\M}_\Ell$. In particular, we have the classical definitions of the $\E_\infty$-rings
\[\Tmf=\O^\top(\overline{\M}_\Ell)\qquad\text{and} \qquad \TMF=\O^\top(\M_\Ell),\]
where $\M_\Ell\subseteq \overline{\M}_\Ell$ is the moduli stack of elliptic curves; see \cite[\textsection7]{ec2name} for an alternative construction of $\O^\top$ on $\M_\Ell$.

\begin{mydef}\label{smfdefinition}
    Define $\Smf=\O^\syn(\overline{\M}_\Ell)$ and $\SMF=\O^\syn(\M_\Ell)$.
\end{mydef}

By \cref{invert_tau_Osyn}, the synthetic $\E_\infty$-rings $\Smf$ and $\SMF$ are synthetic lifts of $\Tmf$ and $\TMF$, and by \cref{signatureofsyntheticglobalsections}, their signature spectral sequences are the DSS for $\Tmf$ and $\TMF$.

\begin{remark}
We will ignore the case of $\tmf$, as there is no clear way to define this object from the perspective used here. The first main problem is that the stack $\M_{\mathrm{cubic}}$ associated to $\tmf$ is not Deligne--Mumford, as the automorphism group of the cuspidal cubic given by $y^2=x^3$ is not an étale group scheme. The second problem is that the map $\M_{\mathrm{cubic}} \to \M_\FG$ is not flat. To see this, we notice that if it were flat, then $\MU_\ast \tmf$ would be flat over $\pi_\ast \MU$, which we can see fails once we know that
\[\MU_\ast\tmf \cong \Z[a_1,a_2,a_3,a_4,a_6, e_n \mid n\geq 4];\]
see \cite[Pr.20.4]{rezktmfnotes}.
One can of course simply \emph{define} $\smf$ to be $\opnu\tmf$.
This definition is partially justified by the observation that the map $\M_{\mathrm{cubic}}\to \M_\FG$ is affine, so that \cref{affine_implies_connective} suggests that $\opnu\tmf$ should be the preferred synthetic lift of $\tmf$.
The $\MU$-synthetic spectrum $\opnu\tmf$ has already been studied by Gheorge--Isaksen--Krause--Ricka \cite{mmf} under the guise of \emph{motivic modular forms}, reinterpreted as an $\MU$-synthetic statement via \cite[Th.7.34]{syntheticspectra}.
\end{remark}

\begin{remark}    
In a similar vein, $\opnu\ko$ is a preferred $\MU$-synthetic model for $\ko$ as the map $\M_{\mathrm{quad}}\to \M_\FG$ (where $\M_{\mathrm{quad}}$ is the stack associated to $\ko$; see \cite[\textsection9]{tmfbook}) is also affine.
In this case, the first two authors previously showed \cite[Pr.4.25]{syntheticj} by direct computation that there are equivalences
\[
    \opnu\KO \cong (\opnu\KU)^{h C_2} \qquad \text{and} \qquad \opnu \ko \cong \tau_{\geq 0}^\uparrow \opnu\KO,
\]
which further confirm this suspicion that $\opnu\ko$ is the preferred synthetic lift of $\ko$; here $\tau_{\geq }^\uparrow$ is the \emph{vertical $t$-structure} of \cite[Th.3.2]{syntheticj}.
We repeat the warning from \cite[Warn.4.26]{syntheticj} that these equivalences do \textbf{not} hold in $\Syn_{\mathbf{F}_2}$.
\end{remark}

Our formalism also extends to synthetic modular forms with level structures.

\begin{variant}
    It is also possible to alter \cref{smfdefinition} to accommodate moduli stacks of elliptic curves with \emph{level structures}. Indeed, for a congruence subgroup $\Ga\subseteq \SL_2(\Z)$, write $\M_\Gamma$ for the moduli stack of elliptic curves with $\Ga$-level structure.
    The map $\M_\Gamma \to \M_\Ell$ is \'etale, and we define
    \[
        \SMF_\Ga \defeq \O^\syn(\M_\Ga).
    \]
    This is a synthetic $\E_\infty$-ring that lifts $\TMF_\Gamma$ and whose signature implements the DSS for $\TMF_\Gamma$.
    For compactifications of $\M_\Ga$, we need to go to the \emph{log \'{e}tale} site of $\overline{\M}_\Ell$ as in \cite{hilllawson}. This is due to the structure map $\overline{\M}_\Ga\to \overline{\M}_\Ell$ generally not being \'{e}tale, but only log \'{e}tale. In particular, replacing the small \'{e}tale site of a Deligne--Mumford stack $\x$ with its \emph{log \'{e}tale site}, for a chosen log structure on $\x$, all of the results of this article follow \emph{mutatis mutandis}. This leads to a synthetic $\E_\infty$-ring
    \[\Smf_\Ga \defeq \O^\syn(\overline{\M}_\Ga)\]
    that, likewise, lifts $\Tmf_\Ga$ and whose signature induces the DSS for $\Tmf_\Ga$.
\end{variant}

The following is simply \cref{sheafismulocalandtaucomplete}.

\begin{cor}\label{taucompleteandmulocal}
    The synthetic $\E_\infty$-rings $\SMF$ and $\Smf$ are both $\nu\MU$-local {\upshape(}a.k.a.\ hypercomplete{\upshape)} and $\tau$-complete.
\end{cor}

Our earlier results let us understand $\SMF$ in more familiar terms.

\begin{cor}\label{SMFstatements}
    The canonical limit-comparison map of synthetic $\E_\infty$-rings
    \[\opnu\TMF = \nu (\O^\top(\M_\Ell))\xrightarrow{\cong} \O^\syn(\M_\Ell) = \SMF\]
    is an equivalence.
\end{cor}

\begin{proof}
    The map $\M_\Ell\to \M_\FG$ defined by the formal group associated to the universal elliptic curve over $\M_\Ell$ is affine: see the first paragraph in the proof of \cite[Th.7.2]{akhilandlennart}. The above equivalence then follows from \cref{affine_implies_connective}.
\end{proof}

In other words, $\nu$ preserves the limit $\O^\top(\M_\Ell)$ defining $\TMF$. While the synthetic spectrum $\SMF$ is therefore not new, we will continue to denote it by $\SMF$ rather than $\opnu\TMF$, in order to emphasise its construction and the fact that its signature implments the DSS for $\TMF$.

\begin{remark}
The map of synthetic $\E_\infty$-rings $\opnu\TMF_\Ga\xrightarrow{\cong} \SMF_\Ga$ is also an equivalence, for any congruence subgroup $\Ga$. Indeed, the map of stacks $\M_\Ga \to \M_{\Ell,\Z[\frac{1}{N}]}$ is finite \'{e}tale, and hence affine, so \cref{affine_implies_connective} applies again.
\end{remark}

Although the affine map $\M_\Ell\to\M_\FG$ extends over the compactification $\overline{\M}_\Ell$, the resulting map is \textbf{not} affine, but merely quasi-affine. We can then use \cref{theorem:when_Osyn_connective} to determine how far away the map $\opnu\Tmf\to \Smf$ is from being an equivalence, or in other words, how far away the ANSS for $\Tmf$ is from its DSS.

\begin{cor}\label{Smfstatements}
    The synthetic $\E_\infty$-ring $\Smf$ is $(-1)$-connective but not connective.
    In particular, there is a cofibre sequence of synthetic spectra
    \[\opnu\Tmf\to \Smf\to \Sigma^{-1}(\pi_{-1}^\heartsuit \Smf).\]
    Moreover, viewing $\pi_{-1}^\heartsuit (\Smf)$ as an element of $\Comod_{\MU_*\MU}$, we have an isomorphism of graded abelian groups
    \[
        \pi_{-1}^\heartsuit(\Smf) \cong \uH^1(\overline{\M}_\Ell \times_{\M_\FG} \Spec L, \, \omega^{*/2}).
    \]
\end{cor}

In particular, $\Smf$ is a synthetic lift of $\Tmf$ different from $\opnu\Tmf$.

\begin{proof}
    In \cite[Pr.5.1]{tmfhomology}, Mathew computes the DSS for $\MU \otimes \O^\top_\x$ for $\x = \overline{\M}_\Ell$, finding that it is concentrated in filtrations $0$ and $1$.
    Therefore \cref{thm:SDSS} shows that $\nu\MU_{*,1}(\Smf)$ is nonzero while $\nu\MU_{*,s}(\Smf)$ vanishes for all $s\geq 2$.
    More specifically, it shows that we have an isomorphism
    \[
        \nu\MU_{*,1}(\Smf) \cong \uH^1(\overline{\M}_\Ell \times_{\M_\FG} \Spec L, \, \omega^{(*+1)/2}).
    \]
    In particular, $\Smf$ is $(-1)$-connective but not connective.\\
    
    By \cref{invert_tau_Osyn}, the cofibre of the map $\opnu\Tmf \to \Smf$ is the $(-1)$-truncation $\tau_{\leq -1}\Smf$.
    Since $\Smf$ is $(-1)$-connective, this cofibre is therefore the desuspension of the discrete object $\pi_{-1}^\heartsuit(\Smf)$, considered as an object of $\Syn_\MU^\heartsuit \cong \Comod_{\MU_*\MU}$.
    The last claim follows from \cite[Th.4.18]{syntheticspectra}, which says that if $S$ is a synthetic spectrum, then we have an isomorphism of graded abelian groups
    \[
        \nu\MU_{*,s}(S) \cong \pi_{-s}^\heartsuit(S)[s],
    \]
    where $[s]$ denotes the $s$-fold grading shift.
\end{proof}

Classically, we have an equivalence of $\E_\infty$-rings $\Tmf[\Delta^{-24}] \cong \TMF$.
Proving a synthetic version of such a statement is more subtle, because the synthetic global sections functor does not preserve all colimits (see \cref{rmk:no_synthetic_affineness,tmfcounterexample}).
Nevertheless, we can show that it does preserve the specific filtered colimits used to invert elements of filtration zero on the global sections of $\O^\syn$.
This in fact holds for more general even-periodic refinements in the place of $\overline{\M}_\Ell$.

\begin{prop}\label{generalinvertingelementsprop}
    Let $(\X,\O^\top_\X)$ be a tame even-periodic refinement.
    \begin{numberenum}
        \item \label{item:modtau_inv_filtration_zero} Let $f \in \uH^0(\X,\omega^{\otimes d}_\X)$ be an element, and write $D(f)$ for its nonvanishing locus on $\X$.
        Regarding $f$ as an element of $\oppi_{2d,\,0}(\O^\syn_\X(\X)/\tau)$, the natural map
        \[
            \left(\O^\syn_\X(\X)/\tau\right)[ f^{-1} ] \xrightarrow{\cong} \O^\syn_\X(D(f))/\tau
        \]
        is an isomorphism.
        \item \label{item:integral_inv_filtration_zero} Let $g\in \oppi_{2d,\, 0}\O^\syn_\X(\X)$ be an element, and write $\overline{g}\in \uH^0(\X,\omega^{\otimes d}_\X)$ for its image under the edge homomorphism.
        Suppose that the underlying stack $\X$ is noetherian.
        Then the natural map
        \[\O^\syn_\X(\X)[g^{-1}] \xrightarrow{\cong} \O^\syn_\X(D(\overline{g}))\]
        is an isomorphism.
    \end{numberenum}
\end{prop}

\begin{proof}
    For \ref{item:modtau_inv_filtration_zero}, we check that the map induces an isomorphism on bigraded homotopy groups.
    As filtered colimits commute with bigraded homotopy groups, we calculate (using \cref{signatureofsyntheticglobalsections}) the left-hand side to be the $\uE_2$-page of the DSS for $\O^\top_\X(\X)$ with $f$ inverted.
    This agrees with the $\uE_2$-page of the DSS for $\O^\top_\X$ at $D(f)$.
    Indeed, sheaf cohomology commutes with filtered colimits in the abelian $1$-category of quasi-coherent sheaves; see \cite[Tag~\href{https://stacks.math.columbia.edu/tag/0GQV}{0GQV}]{stacks}.
    As $f$ is of filtration zero, it comes from a map in the abelian $1$-category of quasi-coherent sheaves.
    This $\uE_2$-page in turn agrees with the homotopy groups of $C\tau \otimes \O^\syn_\X(D(f))$, and we see that the map induces the natural isomorphism on bigraded homotopy groups.

    For \ref{item:integral_inv_filtration_zero}, it suffices to check that the map is an isomorphism after tensoring with $C\tau$ and after $\tau$-inversion.
    The former of the two follows from \ref{item:modtau_inv_filtration_zero} for $f = \overline{g}$, because tensoring with $C\tau$ preserves colimits.
    The latter of the two follows from from \cite[Th.4.14]{akhilandlennart}; note that this result requires $\X$ to be noetherian and the map $\X \to \M_\FG$ to be tame.
    Indeed, this result says that the global sections functor
    \[
        \Ga(\X,-)\colon \QCoh(\X,\O_\X^\top)\to \Mod_{\O^\top_\X(\X)}(\Sp)
    \]
    preserves colimits; in particular, we obtain the chain of natural isomorphisms
    \[
        \Ga(\X,\, \O^\top_\X)[f^{-1}] \cong \Ga(\X,\, \O^\top_\X[\overline{f}{}^{-1}]) \cong \Ga(D(\overline{f}),\, \O^\top_\X).\qedhere
    \]
\end{proof}

Let ${\Delta}^{24} \in \oppi_{576,0} \Smf$ denote the unique element that $\tau$-inverts to the usual periodicity element $\Delta^{24}\in \oppi_{576}\Tmf$; see \cite[Lm.8.1]{smfcomputation} for this computation. 
Applying the above to $\Delta^{24}$, we obtain a synthetic lift of the usual relationship between $\Tmf$ and $\TMF$.

\begin{cor}\label{invertingdelta}
    The natural map of synthetic $\E_\infty$-$\Smf$-algebras $\Smf[{\Delta}^{-24}]\xrightarrow{\cong} \SMF$ is an equivalence.
\end{cor}

\begin{proof}
    This follows from \cref{generalinvertingelementsprop} and the fact that $\M_\Ell$ is exactly the nonvanishing locus of $\overline{\Delta}^{24}$ inside $\overline{\M}_\Ell$. 
\end{proof}

Note that the proof of this statement does not depend on \cref{SMFstatements}, but only uses that $\SMF$ is a synthetic lift of $\TMF$, and that its signature implments the DSS for $\TMF$.\\

We can also compare \cref{smfdefinition} with another possible definition of $\Smf$, one originally proposed by the second author.
Recall that $\tmf=\tau_{\geq 0}\Tmf$ and that $c_4\in \oppi_8\tmf$ is the uniquely defined $\overline{\kappa}$-torsion class detecting the normalised Eisenstein series of weight~$4$.

\begin{prop}
    There is a pullback diagram of synthetic $\E_\infty$-rings
    \[\begin{tikzcd}
        {\Smf}\ar[d]\ar[r]  \pullback &   {\opnu\TMF}\ar[d] \\
        {\opnu\tmf[{c}_4^{-1}]}\ar[r]  &   {\opnu\tmf[{c}_4^{-1}, {\Delta}^{-24}].}
    \end{tikzcd}\]
\end{prop}

\begin{proof}
    The following diagram of synthetic $\E_\infty$-rings
    \[\begin{tikzcd}
    {\Smf}\ar[r]\ar[d]              \pullback &   {\Smf[{\Delta}^{-24}]}\ar[d]    \\
    {\Smf[{c}_4^{-1}]}\ar[r]    &   {\Smf[{c}_4^{-1},{\Delta}^{-24}]}
    \end{tikzcd}\]
    is a pullback. Indeed, by \Cref{generalinvertingelementsprop}, we can identify this square with $\O^\syn(D(f))$ for $f=1,\Delta, c_4$, and $c_4\Delta$, respectively, reading from left-to-right then top-to-bottom. The fact that $\overline{\M}_\Ell = D(\Delta)\cup D(c_4)$ and the fact that $\O^\syn$ is a Zariski sheaf then yields the above pullback of synthetic $\E_\infty$-rings. We are left to identify this square with the square in the statement of the proposition.
    By \cref{invertingdelta} and \cref{SMFstatements}, we see that $\Smf[{\Delta}^{-24}]\cong \SMF \cong \opnu\TMF$, respectively. It suffices to now show that the composite morphism of synthetic $\E_\infty$-rings
    \[\opnu\tmf\to \opnu\Tmf\to \Smf\]
    becomes an equivalence after inverting ${c}_4$. After inverting $\tau$, this follows from the fact that $\tmf=\tau_{\geq 0}\Tmf$ by definition.
    After killing $\tau$, we note that by \cite[Cor.5.3]{tmfhomology} and \cref{signatureofsyntheticglobalsections}, we can identify the bigraded homotopy groups of $\opnu\tmf\otimes C\tau$ and $\Smf\otimes C\tau$ as the $\uE_2$-page of a DSS for $\M_{\mathrm{cub}}$ and $\overline{\M}_\Ell$, respectively.
    As $\overline{\M}_\Ell$ can be written inside $\M_{\mathrm{cub}}$ as the union of the open substacks $D(c_4)\cup D(\Delta)$, and filtered colimits commute with cohomology (see \cite[\href{https://stacks.math.columbia.edu/tag/0GQV}{0GQV}]{stacks}), we obtain the desired equivalence after tensoring with $C\tau$.
\end{proof}


\addcontentsline{toc}{section}{References}
\scriptsize
\bibliography{references}
\bibliographystyle{alpha}

\end{document}